\documentclass[11pt]{article}
\usepackage{amsfonts}
\usepackage{fullpage}
\usepackage{amssymb}
\date{}
\def\C{{\bf C}}

\def\NN{{\mathbb N}}
\def\HH{{\mathbb H}} 
\def\RR{{\mathbb R}} 
\def\TT{{\mathbb T}} 
\def\ZZ{{\mathbb Z}} 
 
\def\CC{{\mathbb C}} 
 \null 

\def\Chat{\hat\CC}
\def\TTT{\cal T}

\def\FFF{{\cal F}}

\def\PPP{{\cal P}}

\def\HHH{{\cal H}}
\def\AAA{{\cal A}}
\def\BBB{{\cal B}}
\def\CCC{{\cal C}}

\newtheorem{thm}{Theorem}[section]
\newtheorem{lemma}[thm]{Lemma}
\newtheorem{cor}[thm]{Corollary}
\newtheorem{prop}[thm]{Proposition}
\newtheorem{remark}{Remark}[section]
\newtheorem{conj}{Conjecture}

\newenvironment{proof}{{\sc Proof.}}{$\;\square$ \vskip .2in}

\def\Definition{\smallbreak\noindent{\bf Definition: \ }}
\def\Definitions{\smallbreak\noindent{\bf Definitions: \ }}
\def\Notation{\smallbreak\noindent{\bf Notation: \ }}

\def\f#1{{f^{\circ #1}}}
\def\fl#1#2{{f^{\circ #1}_{#2}}}
\input{psfig}
\def\IMSmarkvadjust{0 pt}
\def\IMSmarkhadjust{0 pt}
\def\IMSmarkhpadding{0 pt}
\def\IMSpubltext{Published in modified form:}
\def\SBIMSMark#1#2#3{
 \font\SBF=cmss10 at 10 true pt
 \font\SBI=cmssi10 at 10 true pt
 \setbox0=\hbox{\SBF \hbox to \IMSmarkhpadding{\relax}
                Stony Brook IMS Preprint \##1}
 \setbox2=\hbox to \wd0{\hfil \SBI #2}
 \setbox4=\hbox to \wd0{\hfil \SBI #3}
 \setbox6=\hbox to \wd0{\hss
             \vbox{\hsize=\wd0 \parskip=0pt \baselineskip=10 true pt
                   \copy0 \break%
                   \copy2 \break%
                   \copy4 \break}}
 \dimen0=\ht6   \advance\dimen0 by \vsize \advance\dimen0 by 8 true pt
                \advance\dimen0 by -\pagetotal
	        \advance\dimen0 by \IMSmarkvadjust
 \dimen2=\hsize \advance\dimen2 by .25 true in
	        \advance\dimen2 by \IMSmarkhadjust

  \openin2=publishd.tex
  \ifeof2\setbox0=\hbox to 0pt{}
  \else 
     \setbox0=\hbox to 3.1 true in{
                \vbox to \ht6{\hsize=3 true in \parskip=0pt  \noindent  
                {\SBI \IMSpubltext}\hfil\break
                {\it  Conf. Geom. \& Dynam.}~{\bf 1} (1997), 28--57
 
                \vfill}}
  \fi
  \closein2
  \ht0=0pt \dp0=0pt
 \ht6=0pt \dp6=0pt
 \setbox8=\vbox to \dimen0{\vfill \hbox to \dimen2{\copy0 \hss \copy6}}
 \ht8=0pt \dp8=0pt \wd8=0pt
 \copy8
 \message{*** Stony Brook IMS Preprint #1, #2. #3 ***}
}

\begin{document}

\title{Dynamics of the family $\lambda \tan z$ }
\author{Linda Keen\\Mathematics Department\\CUNY Lehman College\\Bronx,
   NY 10468, U.S.A.
\thanks{Supported in part by NSF GRANT DMS-9205433, PSC-CUNY Award,
        I.B.M. and M.S.R.I.} \and
Janina Kotus\\Institute of Mathematics
    \\Technical University of Warsaw\\00-661 Warsaw, Poland
\thanks {Supported in part   by Polish KBN Grants No
21046910 ``Iteracje i fraktale'' and  KBN Grant No 2 P03A  025 12 
"Iteracje funkcji holomorficznych" , the Fulbright Foundation, I.M.S.
Stony Brook  and M.S.R.I.}
}

\maketitle
\SBIMSMark{1995/8}{May 1995}{Revised version: June 1997}

\begin{abstract}
We study the  the tangent family
$\FFF = \{\lambda \tan z, \lambda \in \CC - \{0\}\}$ and 
 give a complete classification of their stable behavior. We also
characterize the the hyperbolic components and 
 give a combinatorial description their deployment in the
parameter plane. 
\end{abstract}

\section{Introduction}

One of the central questions in conformal dynamics is characterizing, in
``natural'' families of meromorphic functions such as $\{z^2 +
\lambda, \lambda \in \CC\}$, those members that define hyperbolic
dynamical
systems.  The Mandelbrot set shows how the hyperbolic systems fall into
connected components of the parameter plane.  Systems within a given
component are topologically conjugate and the  combinatorial structure
of these components has been studied in detail by many people.    McMullen's
survey, \cite{McCan}, and the references therein provide a good
introduction to this work.  

Transcendental meromorphic functions, and $\tan z$ in particular, offer
another interesting one complex dimension domain to study.  
The dynamical properties of the tangent family were first studied by
Devaney and Keen in \cite{DK1,DK2} and pursued by Baker, Kotus, and L\"u
in \cite{BKL2,BKL1,BKL3,BKL4}.  W.H. Jiang  did the computer studies of the
parameter plane in figures 1,2 and 3 as part of his unpublished thesis
\cite{Jiang}.  
These figures were the motivation for our work.

In the quadratic family, one can show that each hyperbolic component in the
parameter plane contains a unique point (its ``center'') that has
a superattractive periodic cycle.  One can enumerate these centers in terms
of  a combinatorial description
of their dynamical behavior.    
Since $\lambda\tan z$ has  no critical
points, and so no superattractive cycles, we need another way
to enumerate its hyperbolic components.  We use parameter values for which
the asymptotic values eventually land on a pole.  We call these parameter
values {\em virtual centers}, which is an apt term dynamically, but a
little dicey geometrically, because they turn out to lie on the boundary of
the hyperbolic components.
We shall see that for each $\lambda$, there is a natural map from the
prepoles of $\lambda \tan z$ to the set of virtual centers.  We use this
correspondence to characterize the hyperbolic components.

The paper is organized as follows.  In section 2, we give the basic
facts about iteration of meromorphic functions.  In section 3, we discuss 
the concepts of   
 quasiconformal
stability  and J-stability of meromorphic functions in terms of
 extensions of holomorphic motions \cite{MSS,McS}.
  In section
4, we  give the special properties 
of the tangent family $\FFF$ and in section 5, we classify the stable
behavior of its 
functions.  In addition, we characterize the connectivity properties of the
Julia sets of functions in $\FFF$ for various values of the parameter
$\lambda$. 
In section 6, we give an alternative characterization for the Julia set of
the tangent functions as the closure of the set of prepoles.  We describe
the combinatorial structure of the prepoles and of the periodic points. 
In section 7 
we  prove that for the tangent family, the set of J-stable parameters
coincides with the set 
of quasiconformally stable parameters.    In
sections 8 and 9 
we describe the deployment of the hyperbolic components and establish the
correspondence between the prepoles and virtual centers.  We exploit this
correspondence, together with the results in section 6, to transfer
properties from the dynamic plane to the parameter plane to obtain the
combinatorial picture of the hyperbolic components.

The authors wish to thank Anna Zdunik for reading preliminary versions of
this paper and for her helpful comments.  They would also like to thank Scott
Sutherland for his help with the figures.  In addition, they thank the
referee for his suggestions for improving the exposition and for finding
errors in our proof of the density conjecture in  earlier versions of this
paper.   The second
author would 
like to acknowledge the hospitality of the Graduate Center at CUNY and the
Institute of Mathematical Sciences at SUNY Stony Brook.  Finally, both authors
thank the Mathematical Sciences Research Institute for its support.

\section{Basics}

\subsection{Meromorphic functions}

 If $f$ is a meromorphic function, the orbits of points fall into three
categories: they may be infinite, they may become periodic and hence
consist of a finite number of distinct points, or they may terminate at a
pole of $f$. For transcendental meromorphic functions with more than one
pole, it follows from Picard's theorem that the set of prepoles, that is,
${\cal P} = \cup_{n=0}^{\infty} \f{-n}(\infty)$ is infinite.
To study the dynamics, we define the {\em stable set}, or {\em Fatou set},
$F_f$ as the set of those points $z$ such that the sequence $\f{n}(z)$
 is defined and meromorphic for all $n$, and forms a normal family in a
neighborhood of $z$.  The  {\em unstable set}, or {\em Julia set}, $J_f$
 is the complement of the Fatou set. Thus $F_f$ is open, $J_f$
is closed and $\bar{\cal P}
\subset J_f$ so that $J_f \neq \emptyset$.
 It is easy to see that both
$F_f$ and $J_f$ are completely invariant.

Meromorphic functions may have values $z \in \hat\CC$ with finite backward
orbit. These constitute the {\em exceptional set} $E_f$. By Picard's
theorem, $|E_f| \leq 2$ and $\infty \not\in E_f$. As for rational
functions, for any $z\in \hat\CC -  E_f$,
 $J_f \subset \overline{\cup_{n=0}^{\infty}\f{-n}(z)}$ and if
 $z \in J_f$, $J_f = \overline{\cup_{n=0}^{\infty}\f{-n}(z)}$. It follows that
 $J_f = \overline{\PPP}$, a fact which is essential to our proofs. We also
 may characterize $J_f$, as we do for rational maps, as the closure of
 the repelling periodic points  \cite{DK1,BKL1}.

The {\em singular set} $S_f$ of a meromorphic function $f$ consists of
those values at which $f$ is not a regular covering.  Therefore at a
singular value $v$ there is a branch of the inverse which is not
holomorphic but has an algebraic or transcendental singularity. If the
singularity is algebraic $v$ is a critical value whereas if it is
transcendental, there is a path $\alpha:[0, \infty) \rightarrow \C$ such
that $\lim_{t\to \infty} |\alpha(t)| = \infty$ and $\lim_{t \to \infty}
f(\alpha(t)) = v$; in this case $v$ is called an {\em asymptotic value} for
$f$. If we can associate to a given asymptotic value $v$ an {\em asymptotic
tract}, that is, a simply connected unbounded domain $A$ such that $f(A)$
is a punctured neighborhood of $v$, then $v$ is called a {\em logarithmic
singularity}.  Note that logarithmic singularities and critical values may
belong to the Fatou set. The functions in the tangent family $\FFF$ have no
critical points and so the singular set consists of asymptotic values
which, in this case, are the points $\pm \lambda i$ of the exceptional set.
These are logarithmic singularities with asymptotic tracts contained in the
upper and lower half planes respectively.

 Let $D$ be a component of the Fatou set; $f$  maps $D$ to a
component, but if the image contains an asymptotic value, the map
may not be onto.  In any case, we call the image component $f(D)$ and note
 the dichotomy:
\begin{itemize}
\item there exist integers $m \neq n   >0$ such that
$\f{n}(D)
= \f{m}(D)$, and $D$ is called {\em eventually periodic},
or
 \item for all $m \neq n$, $\f{n}(D) \cap \f{m}(D) = \emptyset$,
 and $D$ is called {\em a wandering domain}.
\end{itemize}
 The qualitative and quantitative description of the eventually
periodic stable behavior of meromorphic functions is slightly more
complicated than  that of  rational  maps 
 because of the essential singularity at $\infty$
   and the possibility that $f^p$ may not be defined at some values.
Suppose $z_0, z_1=f(z_0), \ldots ,
z_p = f(z_{p-1})=z_0$ is
a periodic cycle for $f$.  The {\em eigenvalue} or {\em multiplier}
 of the cycle is defined to be $\alpha =(\f{p})^{\prime}(z_i)$,
 $i=0, \ldots, p-1$.

 Suppose now that the domain $D$ lands on
 a periodic cycle of domains $D_0, D_1,\ldots D_{p-1}$;
 then either
\begin{enumerate}
 \item $D_i$ is {\em attractive}: each $D_i$ contains a point
 of a periodic cycle with eigenvalue $|\alpha|<1$ and all points
in the $D_i$ are attracted to this cycle. Some domain in this cycle
must contain a critical or an asymptotic value. If $|\alpha|=0$,
the critical point itself belongs to the periodic cycle and the
 domain is called superattractive.
\item $D_i$ is {\em parabolic}: the boundary of each $D_i$ contains
 a point of a periodic cycle with eigenvalue
 $ \alpha = e^{2 \pi i p/q}$,
$(p,q)=1$, and all points in the domains
 $D_i$ are attracted to this cycle. Some domain in this cycle must
 contain a critical  or an asymptotic value.
\item $D_i$ is a {\em Siegel disk}: each $D_i$ is contains a point
 of a periodic cycle with eigenvalue $\alpha=e^{2 \pi i \theta}$
 where $\theta$ is irrational. There is  a holomorphic homeomorphism
mapping each $
D_i$  to the unit disk $ \Delta$, and conjugating
 the first return map $\f{p}$ on $D_i$ to an irrational rotation of
  $\Delta$. The preimages under this conjugacy  of the circles
$|\zeta|=r, \, r<1$ foliate the disks $D_i$ with $\f{p}$ forward invariant
leaves on which $\f{p}$ is injective.
\item $D_i$ is a {\em Herman ring}: each $D_i$ is holomorphically
homeomorphic to a standard annulus and the first return map is
 conjugate to an irrational rotation of the annulus by the
holomorphic homeomorphism. The preimages under this conjugacy
 of the circles $|\zeta|=r,\, 1<r<R$ foliate the disks with $\f{p}$
 forward invariant leaves on which $\f{p}$ is injective.
\item $D_i$ is an {\em essentially parabolic  or Baker} domain:
 the boundary of each  $D_i$ contains a point
 $ z_i$ (possibly $\infty$), $\f{np}(z)\to z_i$  for all
$z \in D_i$, but $\f{p}$ is not holomorphic at $z_i$. If $p=1$, then
  $ z_0=\infty$ .

For examples of essentially parabolic domains with $p\neq 1$
 and $z_i\neq \infty$  see \cite{BKL2}.
\end{enumerate}

\Definition  A meromorphic map $f$ is called {\em hyperbolic} if it is
expanding on its Julia set; that is, there exist constants $c > 0$ and 
 $K > 1$ such that for all $z$ in a neighborhood $V \supset J$,
$|(\f{n})'(z)|>cK^n$. 

\section{Quasiconformal Conjugacy and J-stability}

An effective method for studying the analytic structure of the parameter
space of a family of rational maps or meromorphic functions with finite
singular set is to use the theory of holomorphic motions.  This is
presented in detail for rational maps in \cite{McS}. 
We indicate here how to adapt it for meromorphic functions with finite
singular set.  This is not substantially harder than restricting ourselves
to the tangent family $\FFF$ and gives a more general result.

\bigskip
\Definitions
\begin{itemize}
\item
 A {\em holomorphic family} of meromorphic maps $f_x(z)$ over a
complex manifold $X$ is a holomorphic map $X \times \Chat \rightarrow
\Chat$, given by $(x, z) \mapsto f_x(z)$. 
\item
 A point $x_0 \in X$ is {\em
topologically stable} if there is a neighborhood $U$ of $x_0$ such that for
all $y \in U$ there is a homeomorphism $\phi:\Chat \to \Chat$ such that
$f_y= \phi^{-1}\circ f_{x_0} \circ \phi$.  Denote by $X^{top} \subset X$
the subset of topologically stable parameters. 
\item
 The subset $X^{qc}\subset
X^{top}$ is defined similarly, except that $\phi$ is required to be
quasiconformal.  
\item
The subset $X^{post} \subset X$ of {\em postsingularly stable} points 
consists of those parameter values
for which the combinatorial properties of the singular values is
locally constant; that is,  in 
a neighborhood of $x \in X^{post}$, any relations among the 
forward orbits of the singular values of $f_{x}$ persist. 
For $x \in X^{post}$
the singular points are locally defined distinct holomorphic 
functions of $x$.  Clearly, $X^{post}$ is open and dense in the
subset of $X$ on which the number of singular values counted
without multiplicity is constant and $X^{top} \subset X^{post}$. 
\item
A  {\em
holomorphic motion} of a set $V \subset \Chat$ over a connected complex
manifold with basepoint 
$(X,x_0)$ is a map $\phi: X \times V \rightarrow \Chat$ given by $(x,v)
\mapsto \phi_x(v)$ such that  
\begin{enumerate}
\item for each $v \in
V$, $\phi_x(v)$ is holomorphic in $x$,  
\item
for each $x \in V$, $\phi_x(v)$ is an injective function of $v \in V$,
and, 
\item  at $x_0$, $\phi_{x_0}=v(x_0)$.
\end{enumerate}
\item
A holomorphic motion over a holomorphic family {\em respects the dynamics}
if $\phi_x (f_{x_0}(v)) = f_x(\phi_{x}(v))$ whenever both $v$ and $f_x(v)$
belong to $V$.  
\item
The set $X^{stab}$ of {\em J-stable} parameters is the subset of $X$ where
the Julia set moves by a holomorphic motion respecting the dynamics. 
\end{itemize}

The following theorem is proved in \cite{MSS} (Theorem B).  We state it
here as it applies in our context:

\begin{thm} In any holomorphic family of meromorphic maps with finite
singular set defined  over the complex manifold $X$, the set of J-stable
parameters coincides with the set on which the total number of attracting
and superattracting cycles of $f_x$ is constant on a neighborhood of $x$. 
\end{thm}

As a corollary we have:
\begin{cor} The J-stable parameters are open and dense in $X$. 
 \end{cor}

\begin{proof} Let $N(x)$ be the number of attracting cycles of $f_x$.  By
Fatou  
\cite{Fatou}  
the set on 
which $N(x)$ is a local maximum is open and dense. By the
above theorem, this set coincides with $X^{stab}$.  
\end{proof}

The following theorem is proved in \cite{McS} for rational maps, but the
proof works just as well for meromorphic maps with finite singular set.  

\begin{thm}
\label{thm:qcdense} In any holomorphic family of meromorphic maps with finite
singular set defined  over the complex manifold $X$, $X^{qc}=X^{top}$ and
$X^{qc}$ is open and dense.
\end{thm}

\begin{proof} 
  Fix $x_0 \in
X^{post}$ and let $Y$ be a neighborhood of $x_0$ on which the
combinatorics of the singular values are constant.  Define a
holomorphic mapping $\phi: Y \times S_{f_{x_0}} \rightarrow \C$ by
$(x,v) \mapsto \phi_x(v) = v(x)$.  Note that $\phi$ is a {\em
holomorphic motion} because the singular values move injectively. 

Since the combinatorics of the singular values are constant in $Y$,
the extension of the motion to the forward orbits of the singular
values by $$\phi_x(f_{x_0}(v)) = f_x ( \phi_x(v))$$ is well defined,
injective and depends holomorphically on $x$.  Next, we can extend
this motion to the grand orbits\footnote{Two points, $z,w$ are in the
same {\em grand orbit} of $f$ if there are integers $m,n >0$ such
that $f^{\circ m}(z) = f^{\circ n}(w)$.} of the singular values.
Again, the persistence of the orbit relations is precisely what is 
needed to do this 
 so that the appropriate inverse branch can be
uniquely chosen by analytic continuation.  
Indeed, suppose that  $FS_{f_{x}}$ is the post-singular set and
that $\phi_x(z)$ is defined for some $z \in FS_{f_{x}}$ and $x \in
Y$.  To show that the inverse branches are holomorphic single-valued
functions defining a holomorphic function of $f_{x}^{-1}(\phi_x(z))$
we consider the graph
$$ G=\{(x,w) | f_x(w)=\phi_x(z) \}$$ defined by the preimages and
show that the cardinality of the preimage is constant
on each  component.   

If $z$ is a critical value, its critical preimages $c_{i_1}, \ldots
c_{i_k}$ have constant multiplicity in $Y$ and each defines a
component of $G$ with cardinality given by the multiplicity.
Non-critical preimages map injectively and define components whose
cardinality is one.  If $z$ is an asymptotic value, the number of
asymptotic tracts is constant in $Y$; these ``omitted'' preimages do
not contribute components to $G$, but the number of such ``missing components''
is constant.  Consequently,  the cardinality of
the preimage is constant on components.  

Now we are ready to extend $\phi_x$ to all of $\Chat$, albeit on
a  neighborhood $U \subset Y$, one-third the size of $Y$, so that 
for each $x \in U$,  $\phi_x$ is  
 quasiconformal and respects the dynamics.   The 
  Bers-Royden harmonic $\lambda$-lemma \cite{BersRoyden} asserts
the existence of 
a unique quasiconformal extension, $\phi:U \times \Chat \rightarrow \Chat$
with the  property that  its  
 Beltrami differential for each $x \in U$, $\phi_{\bar{z}}/\phi_z$, is
harmonic\footnote{A 
Beltrami differential defined on 
$U$ is {\em harmonic} if it can be expressed in local coordinates as
$\mu(z) \frac{d\bar{z}}{dz}= \frac{\bar{\Psi}}{\rho^2}$, where $\Psi dz^2$ is
a holomorphic quadratic differential on $U$ and $\rho|dz|$ is the
Poincar\'e metric on $U$.}. 
 We claim that the Bers-Royden extension $\phi$
respects the dynamics; that is $\phi_x = f_x^{-1} \circ
\phi_x \circ f_{x_0}$.  Set  $\psi_x(z) = f_x^{-1} \circ
\phi_x \circ f_{x_0}(z)$ where    
 the branch of the inverse is chosen
by continuation so that $\psi_{x_0}(z)=z$.   On the post-singular set, by
construction,  
 $\psi_x(z) = \phi_x(z)$ already.    Since $f_{x_0}$ and $f_x$ are
conformal, by the chain rule, the Beltrami coefficient of $\psi_x$ is  the
pullback of the 
harmonic Beltrami coefficient of $\phi_{x_0}$.  Moreover, since $f_{x_0}$ is a
hyperbolic isometry off the post-singular set,  the 
Poincar\'e metric is invariant under the pullback.  It follows that the
Beltrami differential of $\psi_x$ 
is again harmonic and by the uniqueness of the Bers-Royden extension $\phi_x
= \psi_x$.  Thus $\phi$ respects the dynamics and $x_0 \in X^{qc}$. 

Therefore we have proved, $X^{post} \subset X^{qc}$.  Since by
definition, $X^{qc}
\subset X^{top} \subset X^{post}$, we conclude $X^{post}=X^{top}=X^{qc}$.
  \end{proof}

\section{The family $\FFF$}
\subsection{Topological characterization}

The family $\FFF =\{f_{\lambda}=\lambda\tan z\} $ is a holomorphic family
over the punctured complex 
plane $X=\C - \{0\}$.  The discussion below shows  that it is also
topologically closed in the following sense:  if $g$ is topologically
conjugate to $f_{\lambda} \in \FFF$ then $g$ is affine conjugate to
$f_{\lambda}$.

Functions with constant Schwarzian derivative $2k$ have two
 asymptotic values and no critical points. If $f$ is such a function it has
 the form
$$f(z)= \frac{A e^{2kz} +B}{C e^{2kz} +D},\,\,\, A,B,C,D \in \CC,
 \,\,AD-BC\neq 0 $$
and the asymptotic values are $A/C$ and $B/D$.
This follows from  the following fundamental theorem of Nevanlinna
(\cite{Nev}, chap. XI, see also \cite{DK1,DK2}).

\begin{thm}[Nevanlinna]
 Meromorphic maps  whose Schwarzian  derivative  is a polynomial of degree
$p-2$ 
are exactly those functions  which have $p$ logarithmic singularities
$a_0,a_1,\ldots,a_{p-1}$. The  $a_i$ need  not be distinct. There are
exactly $p$ disjoint sectors $W_0,\ldots,W_{p-1}$ at  $\infty$,
 each with angle $\frac{2 \pi}{p}$, and  a collection of disks $B_i$, one
 around each
$a_i$,  satisfying:\\
\begin{quote} $f^{-1}(B_i -  \{a_i\} )$ contains a unique
unbounded component $U_i$,  contained in $W_i$, called its asymptotic tract,
 and\\ $f:U_i \mapsto
  B_i - \{a_i\}$ is a universal  covering.
\end{quote}
\end{thm}
  We have, as almost immediate corollaries to this theorem (see \cite{DK2}
for proof).

 \begin{cor} If $f$ has constant Schwarzian derivative $2k$
and if its asymptotic values are $0$ and $\infty$ then $f(z)=\lambda
e^{kz}$ for some $\lambda \in \CC$.
\end{cor}

\begin{cor}
\label{thm:prenev}
 If $f$ has constant Schwarzian derivative $2k$ and if its asymptotic
values are symmetric with respect to the origin and $0$ is fixed then
$f(z)= \lambda \tan z$ for some $\lambda \in \CC$ and $f \in \FFF$ and the
asymptotic values are $\pm \lambda i$.
\end{cor}

In addition we have the topological closure of the tangent family,

\begin{cor}
\label{thm:nev}
If $g(z)$ is topologically conjugate to $f(z)=\lambda \tan z$,
 then $g(z)$ is affine conjugate to $\lambda'
\tan z$ for some $\lambda'$.  Moreover, if $\phi(\lambda i) =
-\phi(-\lambda i)$, then $g(z) =\lambda'\tan z$.
\end{cor}

\begin{proof} Let $g = \phi \circ f \circ \phi^{-1}$ and set $J(z) = -z$.  Then
since $f=J^{-1} \circ f \circ J$, $g = \phi \circ J \circ \phi^{-1}\circ g
\circ \phi \circ J^{-1} \circ \phi^{-1}$.   Now by
theorem~\ref{thm:qcdense} we may assume $\phi$ is quasiconformal.
Therefore, $\hat J = \phi \circ J
\circ \phi^{-1}$ is a holomorphic homeomorphism of $\CC$ of order $2$ and so
must 
have the form $\hat J = -z + b$ for some $b \in \CC$.  

Topological conjugation preserves the asymptotic values so $v =
\phi(\lambda i)$ and  
$v'= \phi(-\lambda i)$ are the asymptotic values of $g$.  Since $J$
exchanges the asymptotic values of $f$, 
$\hat J(v)=v'$.  By hypothesis,
$v'=-v$ so that $b=0$, $\hat J(z) = -z$ and $g(0) = -g(0)=0$.  Applying
corollary~\ref{thm:prenev} we are done.  
\end{proof}

\begin{remark}
 If $f$ has constant Schwarzian derivative $2k$ and if its asymptotic
values are equal, then $k=0$ and $f$ is constant; since it fixes the
origin, $f\equiv 0$.
\end{remark}

\subsection{Symmetry}
\label{sec:symmetry}

 The symmetry of the maps in $\FFF$ with respect to $0$ implies that
the stable and unstable sets are respectively symmetric with respect
to the origin. Precisely,
\begin{prop}
\label{thm:symmetry}
For $f_{\lambda} \in \FFF$,
\begin{enumerate}
\item $z \in J_{\lambda}$ if and only if $-z \in J_{\lambda}$,
 hence $z \in F_{\lambda}$ if and only if
$-z \in
F_{\lambda}$, and
\item
$z \in J_{\lambda}$ if and only if $z \in J_{-\lambda}$
 hence $z \in F_{\lambda}$ if and only if
$z \in F_{-\lambda}$.
\end{enumerate}
\end{prop}

\begin{proof} For $z\in \CC$ we have
$$f_{\lambda}(-z) = - f_{\lambda}(z), $$ so for all $k\in \NN $, $$
\fl{k}{\lambda}(-z)=-\fl{k}{\lambda}(z),\, \mbox {and }\,
[\fl{k}{\lambda}(-z)]^{\prime}=[\fl{k}{\lambda}(z)]^{\prime}. $$ It follows
that if $\fl{p}{\lambda}(z)=z$, then $\fl{p}{\lambda}(-z)=-z. $ Moreover if
$z_0, \ldots z_{p-1}$ is a periodic cycle of period $p$ then either $p$ is
even and $-z_0, \ldots, -z_{p-1}$ is a cyclic permutation of the former
cycle, or it is a distinct cycle of period $p$ and $p$ may be either odd or
even. The symmetry also implies that if $z$ is a repelling, attracting,
parabolic, or irrational neutral periodic point (stable or unstable), $-z$
has the same property.  Note that there are never superattracting cycles
because there are no critical points.  Since the Julia set is the closure of
the repelling 
periodic points, we can say that $ z\in J_{\lambda} $ if and only if $
-z\in J_{\lambda}$.  Consequently $ z\in {F}_{\lambda}$ if and only if $-z
\in {F}_{\lambda}$. In particular, if $\lambda i \in {F}_{\lambda}$ then $
-\lambda i \in {F}_{\lambda}$ and the orbits $\{\f{n}_{\lambda}(\lambda
i)\}$ and $\{\f{n}_{\lambda}(-\lambda i)\},\,n \in \NN$ are symmetric; they
either belong to the immediate basin associated to one attracting (or
parabolic) cycle of even period or they belong to two distinct symmetric
attracting (or parabolic) cycles of $z$ and $-z$.  Note that in both cases
the immediate basins are symmetric.

 To prove 2. we note that for $z\in \CC $  and
 $k\in\NN$
$$\fl{k}{-\lambda}(z) =(-1)^k\fl{k}{\lambda}(z)\,\, \mbox{and}
\,\,  f_{\lambda}(-z)=f_{-\lambda} (z). $$ So if $p$ is even
 $$\fl{p}{\lambda}(z)=z\,\, \mbox{implies} \,\, \fl{p}{-\lambda}(z)=z,$$
while if $p$ is odd
 $$\fl{p}{\lambda}(z)=z\,\, \mbox{implies}\,  \fl{2p}{-\lambda}(z)=z.$$
If $p$ is odd, this means that if $z_0, \ldots z_{p-1}$ and
$-z_0, \ldots
-z_{p-1}$ are distinct cycles for $f_{\lambda}$,
they belong to the same cycle for $f_{-\lambda}$ but the cycle has
 double the period.  For $p$ even the same may be true, but there is
 another possibility: $f_{-\lambda}$ may still have two distinct cycles
 but they are now
 $z_0,-z_{1}, z_2,\ldots,z_{p-2},-z_{p-1}$ and
 $-z_0,z_{1},-z_2,\ldots,-z_{p-2},z_{p-1}$.

 Since
$$[\fl{k}{-\lambda}(z)]^{\prime}=(-1)^k[\fl{k}{\lambda}(z)
]^{\prime} \,\,  \mbox{for}\,\,\,  k\in \NN, $$
 if $z$ is a repelling, attracting, parabolic,
 or irrational neutral  periodic point of $f_{\lambda}$, then $z $ is
a periodic point of the same type for $f_{-\lambda}$
(with possibly double or half the period). To complete the proof
we again use the property that the Julia set is the closure of
the repelling periodic points and the classification of periodic
stable behavior.
 \end{proof}

\begin{remark}
 We see that although the map $J(z)=-z$ satisfies $J \circ f_{\lambda} =
f_{\lambda} \circ J$ for all $\lambda$, it does not conjugate the dynamics
of $f_{\lambda}$ and $ f_{-\lambda}$ because distinct symmetric cycles of
period $p$ of $f_{\lambda}$ often (but not always) become single cycles of
period $2p$ of $f_{-\lambda}$ and vice versa. The map $j(z)=\bar{z}$, on
the other hand, conjugates $f_{\lambda}$ to $f_{\bar\lambda}$ and
conjugates the dynamics properly.
\end{remark}

\subsection{Inverse branches}

   The function $\tan z$ is periodic with period $\pi$ and so is
 an infinite to one cover of $\hat\CC - \{ \pm \lambda i\}$. The
origin is a fixed point and the points
$q_n= n\pi$, $n=\pm 1, \pm 2, \ldots $ are also preimages of
 the origin. The poles are
$s_n =
(n+1/2)\pi$, $n=\pm 1, \pm 2, \ldots$; the image of
the real line between any pair of poles is the whole real line.
  The image of the imaginary axis and its translates by $\pi$
  is the vertical segment of the imaginary axis between $i$ and
 $-i$; the image of the vertical lines
$l_n=(n+1/2)\pi + it, t \in \RR,  n \in \ZZ$ is the pair of
 segments of the imaginary axis above and below $\pm i$. We
 denote by $L_n$ the half open vertical strip between
the
lines $l_{n-1}$ and $l_n$ and containing the line $l_{n-1}$.
 From the above, we see that it is divided into four quadrants,
each a preimage of the corresponding quadrant in $\CC$.
 The function $f_{\lambda}= \lambda \tan z$
therefore maps each strip $L_n$ onto $\hat\CC - \{\pm \lambda i\}$
 and takes the quadrants in $L_n$ onto the quadrants of $\CC$
 formed  by multiplying the real and imaginary axes by $\lambda i$.

   The inverse map of $f_{\lambda}$ is given by the following formula
\begin{equation}
f^{-1}_{\lambda}(z)=\frac{1}{2i} \log \left(\frac{\lambda+ i z }
{\lambda - iz }\right).
\end{equation}
   Let $\lambda= a_1+ia_2,z=x+iy $ then
\begin{eqnarray}
\Re f_{\lambda}^{-1}(z)&=&\frac{1}{2}\arctan \left(\frac{2a_1x +2a_2y}
{|\lambda|^2-|z|^2}\right)\\
\Im f_{\lambda}^{-1}(z)&=& -\frac{1}{4}\log
\left(\frac{(|\lambda|^2-|z|^2)^2+(2a_1x+2a_2y)^2}
{[(a_1+y)^2+(a_2-x)^2]^2}\right),
\end{eqnarray}
 where in equation (2) we must specify which branch of the $\arctan$ we
use. We therefore
  denote by $f_{n,\lambda}^{-1} $ (or $f_{n}^{-1}$
if no confusion will result)  the branch of the inverse map whose real
part   is in the strip $L_n$.  For a given $p\in \NN$ we define a branch of
$f^{-p}_{\lambda}$ by
\begin{equation}
 f_{{\bf n}_p, \lambda}^{-p}=f^{-1}_{n_p, \lambda}\circ \cdots \circ
      f^{-1}_{n_2, \lambda}\circ f^{-1}_{n_1, \lambda},
\end{equation}
and set ${\bf n}_p= (n_1, n_2,\ldots, n_p) $.  We call ${\bf n}_p$ the {\em
itinerary} of the map $f_{{\bf n}_p,\lambda}^{-p} $.

\section{The Dynamic Plane: the Fatou set $F_{\lambda}$ }
\label{sec:dynamic plane}

The maps in $\FFF$ have no critical points so 
 there are no superattracting cycles. Moreover,  it follows  from
 \cite{DK2} that functions in $\FFF$, have
no essentially parabolic domains. In this  section  we
shall prove 
that no function in the family $\FFF$ can have a Herman ring. In \cite{DK2}
and in \cite{BKL4}  it is proved that maps in $\FFF$ have no wandering
domains. The stable behavior for functions in $\FFF$ is therefore
 either eventually attractive, parabolic, or lands on a cycle of Siegel
disks.

It therefore follows, just as for rational maps that an equivalent
characterization of hyperbolicity for maps in $\FFF$ is (see e.g.
\cite{Milnors-notes}): 

\begin{prop}  A map in $\FFF$ is 
 hyperbolic if and only if the closure of
its post-singular set, $\overline{\cup_{n \geq 0,v \in S_f} \f{n}(v)}$,
is disjoint from its Julia set. 
\end{prop}

In analogy with the situation for rational maps, 
we shall prove the following dichotomy exists:

\begin{thm}
\label{thm:dichot}
For $f_{\lambda} \in \FFF$  either
\begin{itemize}
\item $|\lambda|<1 $,  the Julia set $J_{\lambda} $ 
is locally a Cantor set and the Fatou set consists of a single infinitely
connected invariant component, or
\item $ |\lambda|\geq 1 $,  the the Julia set  $J_{\lambda} $
   is  connected and  the Fatou set consists of either exactly two or
infinitely many 
simply connected components.  
\end{itemize}
\end{thm}

\medskip

 Let $D_0, D_1 $ be components of the Fatou set  with $f: D_0
 \to D_1$, where for readability, and provided there is no ambiguity, we omit
the subscript $\lambda$. Since  there are no critical points, $f|_{D_0}$ is
 a covering of its image and the degree is either 1 or infinite.
\begin{prop}
 If the degree of $f|_{D_0}$ is infinite then $D_0$ is unbounded.
\end{prop}

\begin{proof} Suppose $D_0$ is bounded and let $\zeta$ be some point in $f(D_0)$.
Let $w\in \partial{D_0}$ be an accumulation point of $P=\{z_n=
f^{-1}_n(\zeta) , \, n\in \bf Z \}$.  Choose a subsequence $z_{n_j}\in P$
such that $z_{n_j}\to w$. Then, since $f$ is analytic and
$f(z_{n_j})=\zeta, \, f(w)=\zeta.$ This, however, contradicts the
invariance of the Fatou set.
 \end{proof}

It follows that if the degree of $f|_{D_0}$ is 1 then $f(D_0)=D_1$, while
if the degree is infinite and if $D_1$ contains an asymptotic value, then $D_0$ contains an asymptotic tract 
and $f$ is a infinite degree covering map from the asymptotic tract onto a
punctured neighborhood of the asymptotic value.

\begin{cor}
Any attractive or parabolic cycle of components contains an
unbounded component.
\end{cor}

\begin{proof} Some component of the periodic cycle must contain an asymptotic
value hence its preimage contains an asymptotic tract.
 \end{proof}

An immediate corollary is
\begin{cor}
Any completely invariant component is unbounded.
\end{cor}

\begin{remark} It was proved in \cite{BKL3} that a transcendental meromorphic
function with finitely many singular values has a maximum of two
completely invariant components.  Moreover, if $F_{\lambda}$ has two
completely invariant components then each is simply connected.  If the
number of components of $F$ were finite, each would be completely
invariant under a high enough iterate and so the number of components of
$F$ is 1,2, or infinity.
\end{remark}

\begin{prop} The connectivity of any component in the immediate basin of
attraction of an attracting cycle is either 1 or infinity. For a parabolic
cycle it is 1.
\end{prop}

\begin{proof}
Let $D_0$ be such a component, and let $z_0$ be an attractive periodic
point in $D_0$ with first return map $f^p$. Let $U$ be a neighborhood of
$z_0$ on which the linearization is defined and set $U_n=f^{-(np)}(U_{n-1})
\cap D_0$. If the connectivity of $D_0 >1,$ there is a $k$ such that the
connectivity of $U_n >1$ for all $n>k.$ Since the degree of $f^p$ is
infinite, the connectivity of $U_n$ is infinity and $D_0$ is
infinitely connected.

Now consider the petals of the immediate basin of a parabolic cycle.  There
are always at least two petals, either because there are two fixed points,
or because there
is a period $2$ cycle, or, because if there is a single fixed point at the
origin it 
must have two petals.  An argument similar to the above shows that each
petal is either simply or infinitely connected.   An infinitely
connected component is necessarily unbounded and intersects a full
neighborhood of infinity and hence contains both asymptotic tracts.  Since 
this is true for each component of the cycle, there is only one
and it is completely invariant.  
\end{proof}

 Maps in $\FFF$  with $|\lambda|<1$  were discussed in \cite{DK1} where
 the following was proved.

\begin{prop}
\label{thm:originattr}
 If $|\lambda|<1$  then $f$  has a single
attracting fixed point, $0$, its Fatou set consists of one
completely invariant component containing both asymptotic values
 and its Julia set $J_{\lambda}$ is a locally a Cantor set.  The map
 $f|J_{\lambda}$ is conjugate to the shift on infinitely many
 symbols.
\end{prop}

The point $0$ is  fixed for all
 $f_{\lambda}\in \FFF$ and $f^{\prime}_{\lambda}(0)=\lambda.$
If $|\lambda| = 1$ and $0$ is  a stable neutral  point,
  the immediate basin of attraction is the $f_{\lambda}$ invariant Siegel
disk $D_0$ which is simply connected. Since all its preimages are also
simply connected, and there are no other stable phenomena, all components
of the Fatou set are simply connected.

If $0$ is an unstable neutral point and the Fatou set is non-empty the
  origin must be a parabolic fixed point. The immediate attractive basin
consists of one or two cycles of simply connected domains containing
 {\em attractive petals} (see e.g. \cite{Milnors-notes}).

Next we have,

\begin{prop}
\label{thm:originrep} Suppose $|\lambda| > 1$ and
let $D_0,\ldots,D_{p-1}$ be the components of the immediate basin of
attraction of an attractive or parabolic periodic cycle. Then
$D_i, i=1, \ldots p-1$ are simply connected.
\end {prop}

\begin{proof}
Suppose first that the cycle of components contains only {\em one} of the
asymptotic values and suppose $p >1$.  Assume (relabeling if necessary)
that $D_0$ is the component of the cycle containing the
 asymptotic tract. Therefore the degree of $f: D_0 \to D_1 -\{\lambda i\}
$ is infinite while the degree of each $\f{p-k} :D_k \to D_0,\,
k=1,\ldots, p-1$
 is $1$ and each is onto.  If $D_0$ is multiply and hence infinitely
connected, a degree argument
 implies that $D_k$ is also infinitely connected. Suppose now that some
$D_i ,\, 0<i\leq p-1$ is infinitely connected. Since  $\f{p-i}:D_i\to D_0
$ and
 $\f{i-k}:D_k\to D_i ,\, 1 \leq k \leq i-1$ are degree one, all the components
 have the same connectivity as $D_i.$

Let $\gamma$ be a smooth non-homotopically trivial curve in $D_0$; it thus
has a prepole in the bounded component of its complement.  For some $k$,
though, ${\gamma}_k=\f{k}(\gamma)$ does have a pole in the bounded
component of its complement and $\gamma_{k+1}\subset D_{k+1 \mathop{mod}
p}=D_j$ has non-zero winding number about the origin.

Now consider the symmetric cycle of domains (containing the other
asymptotic value) and denote the symmetry by  $J(z)=-z$. The curve
$\f{k+1}(J(\gamma)) \subset J(D_j)$ also
has non-zero winding number about the origin; since
$\f{k+1}(J(\gamma))=J(\gamma_{k+1})$, these curves intersect, $D_j \cap
J(D_j) \neq \emptyset$ and the cycles are not distinct.

If there is only one cycle of components $p$ is even and we may label so
that $D_0$ contains the preasymptotic tract of $\lambda i$ and
$D_{p/2}=J(D_0)$ contains the preasymptotic tract of $-\lambda i$. Taking
$\gamma$, $\gamma_{k+1}$ and $D_j$ as above, $J(\gamma)\subset D_{p/2}$
and $D_j \cap D_{j+p/2 \mathop{mod} p} \neq \emptyset$ so must be equal.
But a single component cannot contain two distinct periodic points, so
there is a single attracting fixed point which must be the origin and
$|\lambda|<1$
contradicting the hypothesis.
 \end{proof}

Finally we prove that there are no doubly connected domains.
\begin{prop}
\label{thm:noherman}
 A function $f_{\lambda} \in \FFF$ cannot have a cycle of Herman rings.
\end{prop}

\begin{proof} Suppose $D_0, \ldots D_{p-1}$ were a cycle of Herman rings of period
$p$ for $f$. For each $i, i=0, \ldots p-1$, the first return map $\f{p}:
D_i \rightarrow D_i$ would be conjugate to an irrational rotation and
therefore have degree one. Let $\gamma$ be an $\f{p}$ invariant leaf of
$D_i$. Since $D_i$ is multiply connected, $\gamma$ contains a preimage of a
pole in the bounded component of its complement, $B_{\gamma}$. It follows
that some iterate $\f{n}(\gamma), n \geq 0, $ contains some pole $s_k$ in
$B_{\f{n}(\gamma)}$.
 Assume this is already true for
$\gamma$ and moreover, that $\gamma$ has been chosen to pass very close to
the pole.

  We claim $B_{\gamma}$ also contains $-s_k$.  If not, by symmetry
$-\gamma$ is an invariant leaf of another ring and
$B_{-\gamma}$ contains $-s_k$.  Both $f(\gamma)$ and
$f(-\gamma)$ must have non-zero winding numbers with respect to the origin,
and must intersect ---but they cannot.  Now since $B_{\gamma}$
contains both $s_k$ and  $-s_k$,
$\gamma$ winds around the
origin and intersects some number of strips $L_n$, $n= -k, -(k-1), \ldots,
0, \ldots, k$.   The winding number of $f^2(\gamma)$ therefore is
 at least $2k$ with respect to each of the asymptotic values. Applying
$f$ another $p-2$ times we see that $\f{p}$ cannot be degree one on
$\gamma$ which is a contradiction.
 \end{proof}

Theorem~\ref{thm:dichot} now follows from
propositions~\ref{thm:originattr},~\ref{thm:originrep} and ~\ref{thm:noherman}.

As a corollary to this discussion we see that the hyperbolic maps in $\FFF$
are precisely those that have an attracting periodic cycle.

\section{The Dynamic Plane:  The Julia set $J_{\lambda}$}

\subsection{Combinatorics of the prepoles}
\label{sec:combprepoles}

In what follows we shall always use the notation
$s_n=(n+1/2)\pi$ to denote a pole of $f_{\lambda}$.

We define the {\em pre-poles of order $p$}, $p \in \NN \cup \{0\}  $,
 as the set
$${\cal P}_p = \{ f^{-p}_{{\bf n}_{p}}(\infty):
 {\bf n}_{p}= (n_1,n_2,\ldots,n_p), n_i \in \ZZ \},$$
 where the inverse branches are defined by  (1)-(4).  Then if $v \in
 \PPP_p$, the  ($p-1)$st-iterate, $\f{p-1}(v)=s_{n_1}$ and
$s_{n_1}$ is the pole $s_{n_1}=f^{-1}_{n_1}(\infty)$. It is clear that the
prepoles of order $p$ are in one to one correspondence with the $p$-tuples
of positive integers.  Set
$\PPP = \cup_0^{\infty}\PPP_p$.

Since the point at infinity is an essential singularity
 for $f_{\lambda}$, every value (except the asymptotic values) is taken
infinitely many times
in a neighborhood of infinity. The domains
${\cal A}^+=\{z |\Im z > r\}$ and ${\cal A}^-=\{z |\Im z <- r\} $
for any $r>0$ are the asymptotic tracts of $\pm \lambda i$
respectively and $A^{\pm}_r = f_{\lambda}({\cal A}^{\pm})$
is a  disk punctured at
$\pm \lambda i$. The preimages of
points outside $A^{\pm}_r$ therefore cluster around the real
 axis $\RR$. The real axis  defines what are classically
known as the {\em Julia directions} for the singularity at
infinity. A prepole $v$ of order $p$ is an essential singularity
 of the function $\f{p+1}$. Its Julia directions are given by
 $f^{-p}_{{\bf n}_{p}}=f^{-1}_{n_p}\circ f^{-1}_{n_2} \circ
 \ldots f^{-1}_{n_1}(\RR)$. We may also consider the preimages
$B^{\pm}_v=f^{-p}_{{\bf n}_{p}}({\cal A}^{\pm})$ of the asymptotic
 tracts at the prepole $v$. We call these the
{\em pre-asymptotic tracts} at $v$. They consist of a pair
of smooth disks, tangent at the prepole, and tangent to the Julia
 directions there.

\begin{prop}
\label{thm:poleitn}
The accumulation points of ${\cal P}_p$ belong to
$\cup_{k=0}^{p-1}{\cal P}_k$. For $p \geq 1$ let $v_n\in {\PPP}_p,\, n \in \ZZ$
be a sequence of prepoles with  itinerary ${\bf n}_{p}
=(n_1, n_2,\ldots, n_p)$; 
\begin{itemize}
\item  if the entries $n_i,\, i=1,2,\ldots,p-1 $ are the same for all
    $v_n$ and if $n_p = n$ then the accumulation point of
    $v_n$ belongs to ${\PPP}_0=\{\infty\}$.
\item  if  the entries $n_i, i=2,\ldots,p  $ are the same for all
    $v_n$ and $n_1=n$, then the accumulation point of
    $v_n$  is a prepole $v \in {\PPP}_{p-1} $ with
      itinerary ${\bf n}_{p-1}=(n_2,\ldots, n_p)$.
\end{itemize}
\end{prop}
\begin{proof} Let $v_n, n=1,2,\ldots$ be a sequence of prepoles in
 ${\cal P}_p$ accumulating at a point $v$. If
 $ v \not\in\cup_{k=0}^{p-1}{\cal P}_k$, then
 $\f{p}(v)$ would be defined. Since the  points $v_n$
 are poles of $\f{p}$, any accumulation point must be a
non-removable singularity for $\f{p}$ and thus $\f{p}$ is not defined
 there.

 The itinerary of  $v_n$ is
     ${\bf n}_{p}
   =(n_1, n_2,\ldots, n_p)$. Assume  the first $p-1$
   entries are the same for all $v_n$ and and assume $n_p=n$.
   By the definition of inverse branches  of $f$, we have that
  $v_n=f^{-1}_n(f^{-1}_{n_{p-1}}\circ \ldots
  \circ f^{-1}_{{n_1}}(\infty))$, so $\Re v_{n}\in L_n$.
   Thus the only  accumulation point $v$  of $v_n$  is $\infty$;
   that is,  $v\in { \PPP}_0= \{\infty\}$.

Suppose now that the
   last $p-1$ entries are the same for all $v_n$ and
     $n_1=n$. It follows that
   $v_n=(f^{-1}_{n_p}\circ\ldots \circ f^{-1}_{n_2})
     \circ f^{-1}_n(\infty)$ and  consequently
  $\lim_{|n|\to \infty }v_n=(f^{-1}_{n_p}\circ\ldots
    \circ f^{-1}_{n_2}) \circ (\lim_{|n|\to \infty }
     s_n)=f^{-1}_{n_p}\circ\ldots
    \circ f^{-1}_{n_2}(\infty) =v \in {\PPP}_{p-1} $
    with itinerary ${\bf n}_{p-1}=( n_2,\ldots, n_p) $.
               \end{proof}

As a corollary we see that, although any point in the Julia set
 is an accumulation point of prepoles, the orders of the prepoles
must go to infinity. Precisely,
\begin{cor}
\label{thm:polesacc}
If $z \in J_f - {\cal P}$ then for each integer $p > 0$ there is
 a neighborhood $U$ of $z$ such that
 $$U \cap \left(\cup_{k=0}^{p-1} {\cal P}_k\right) = \emptyset.$$
\end{cor}

\subsection{Combinatorics of the repelling periodic points }

We have characterized the Julia set in two ways: as the closure
 of the repelling periodic points and as the closure of the
 prepoles. We now show how these two characterizations are related. For
 transcendental entire functions, repelling periodic points
of fixed order $p$ may accumulate only at infinity. For
meromorphic functions, however, they  may have other
 accumulation points as we show below  for functions in $\FFF$.

\begin{prop}
\label{thm:proppp} If $v \neq \pm\lambda i$ is a prepole of order
 $p-1 > 0$ for $f=f_{\lambda}$, there exists a sequence of
points $z_k$, $k=1,2,\ldots,$ such that
\begin{enumerate}
\item $\f{p}(z_k)=z_k$
\item $z_k \rightarrow v$ as $k \rightarrow \infty$
\item if
$    m_k=(\f{p})^\prime(z_k)$ is the  multiplier of the
 periodic cycle containing $z_k$, then
 $|m_k| \rightarrow \infty$ as $k \rightarrow \infty$.
\end{enumerate}
\end{prop}

\begin{proof} Let $s_n$ be the pole such that $\f{p-2}(v)=s_n$ and
let $U=B(v,\epsilon)$ be a neighborhood of $v$ with
$\epsilon < |v-(\pm \lambda i)|$. Then $\f{p-2}(U)$ is a
neighborhood of $s_n$ and, replacing $f$ by its principal part,
 we see that  $K=\f{p-1}(U) \subset\subset \{z : |z| > R \}$.
   For each branch of
 the inverse, we obtain an open set with compact closure,
 $U_k=f^{-1}_k(U) $ contained in the strip $L_k$.  Clearly,
all but finitely many of the $U_k$ are contained in $K$.
Fix $k=k_0$ such that $U_{k} \subset K$ for all
 $k \geq k_0$, and let ${\bf n}_{p-2}$ be the itinerary such that
 $ v=  f^{-( p-2)}_{{\bf n}_{p-2}}(s_n)$.  Then
 $$V_k= \fl{-(p-1)}{{\bf
n}_{p-1}}(U_k)=
\fl{-(p-2)}{{\bf n}_{p-2}}\circ f^{-1}_n(U_k)$$ satisfies
$\overline{V_{k}} \subset {\mathop {Int}}(U)$. Therefore
 $\fl{-p}{{\bf
n}_p}=\fl{-(p-1)}{{\bf n}_{p-1}} \circ f^{-1}_k$
 maps $U$ onto $V_k$, and hence into itself, in a one to one
 fashion. By the Schwarz lemma  there is an
 attracting fixed point $z_k$ of $\fl{-p}{{\bf
n}_p}$ in $U$
 for all $k \geq k_0$ proving 1.

\begin{sloppypar}
Thus, for each $k \geq k_0$ there is
 a periodic cycle
 $$ z_{k,0}=z_k,\,z_{k,1}=f(z_{k,0}),\,\ldots\, z_{k,p-1}=
f(z_{k,p-2}), $$ where the points $z_{k,p-1}$ are contained in
 the sets $U_k$.
Hence as $|k| \to \infty$, $\Re z_{k,p-1}
\rightarrow \infty$ while $\Im z_{k,p-1}$
 remains bounded. Moreover
$z_{k,p-2}=f^{-1}_n(z_{k,p-1}) \rightarrow s_n$ so that
 $z_{k,0} \rightarrow v$, proving part 2.
\end{sloppypar}

 Increasing
$k_0$ if necessary, we may assume that
$|z_{k,0} - (\pm \lambda
i)| \geq \epsilon/2$ for all
 $k \geq k_0$. Using the formula
$f^{-1}(z)= \frac{1}{2i} \log ((\lambda + i z)/(\lambda - i z))$
we compute that $\Im z_{k,p-1} = {\rm O}(|\log \epsilon|)$.
Since for $j=0, \ldots p-2$, the points $z_{k,j}$ are contained inside
 $\f{j}(\overline V_k)$,  we see that
 there is some constant $C$ such that $|\Im
z_{k,j}| < C$
for all $k \geq k_0$. The formula
 for the multiplier of the cycle is
 $$m_k=\lambda^p\Pi_{j=0}^{p-1}
\sec^2{z_{k,j}}.$$ The terms in this
product are all bounded from below since the imaginary parts
of the $z_{k,j}$ are bounded by $C$ or $|\log
\epsilon|$.
For notational simplicity we drop the subscript $k$ and the
$\pm$ for $\lambda$ in these estimates. Write
$z_0 \approx \lambda i + \delta$ where $|\delta|>\epsilon/2$; Then $z_{p-2}
\approx s_n$ and
$z_{p-1}
\approx x+\frac{i}{2}|\log\epsilon|$.
 Hence
$$\lambda \tan(z_{p-1}) \approx
\lambda i + \delta\, \mbox{ and }
 \lambda^2\sec^2(z_{p-1} )\approx 2 i \lambda \delta.$$
Also, $\lambda \tan(z_{p-2})=z_{p-1}$ so
 $\sec^2{z_{p-2}}= 1 + (z_{p-1})^2/\lambda^2 $. It follows that
 $m ={\rm O}(\Re z_{p-1})^2$ proving 3.
 \end{proof}

\Notation
Let $v \in {\PPP}_{p-1}$ be a prepole of order $p-1$ with itinerary $ {\bf
n}_{p-1}= (n_1,n_2, \ldots, n_{p-1})$ and let $z_{k,i},\,\, i=0,\ldots,p-1
\,\,, k \in \bf Z $ , be the sequence of repelling periodic points of
order $p$
defined
in proposition ~\ref{thm:proppp}. We define the {\em itinerary} of $z_{k,0}$
as ${ \bf n}_{p}= (n_1,n_2, \ldots, n_{p-1,k})$
where
\begin{equation}
 z_{k,0} \in L_{n_{p-1}},\, z_{k,1} \in L_{n_{p-2}}, \ldots,
     z_{k,p-2}\in L_{n_1},\, z_{k,p-1}\in L_{n_p}=L_{k}.
\end{equation}
For consistency with our previous definitions for prepoles, the entries
mark the branches of the inverse used to traverse the cycle backwards.
Similarly, for any periodic cycle, we define its  itinerary in terms
of the  
branches of the inverse used to traverse the cycle backwards.

\subsection{The Julia set for special $\lambda$'s}

We finish our discussion of the dynamic plane with

\begin{prop} \label{thm:prop2.11}
Assume that $f_{\lambda} \in \FFF$.
\begin{enumerate}
\item If the asymptotic value $\lambda i \in \PPP$  then $J_f=\hat\CC$.
\item If the asymptotic value $\lambda i \in \PPP_{p-1}$ then
there is a neighborhood $U$ of $\lambda i$ that does not contain
 repelling points of period $p$.
\end{enumerate}
\end{prop}
\begin{proof} By symmetry, if $\f{p-1}(\lambda i)=s_n$ for some pole
 $s_n$, then  $\f{p-1}(-\lambda i)=-s_n$, so $-\lambda i$ is also
a prepole of the same order as $\lambda i$. Consequently $J_f=\hat\CC$,
since by the classification theorem, the asymptotic values must have
infinite forward orbits for any stable behavior to occur.

 The proof
 of 2. follows from the construction in the proof of proposition
 \ref{thm:proppp} of  repelling
periodic points accumulating at a prepole  that is not an asymptotic
value.
 \end{proof}

We state the  following theorems here for functions in $\FFF$ although
they hold for more general classes of meromorphic functions.
Their proofs can be found in \cite{DK1} and \cite{KK} respectively.

\begin{thm}[\cite{DK1}]
If $\lambda i \in \PPP$ then $J_{\lambda}$ contains a forward invariant set
called a {\em Cantor bouquet}(defined in \cite{DevKr}).
\end{thm}

\begin{thm}[\cite{KK}]
 If $\lambda i$ belongs to a repelling periodic cycle then
$f_{\lambda}$ acts ergodically on $J_{\lambda}$.
\end{thm}

\section{The parameter space: Analytic structure}

\subsection{J-stability and quasiconformal conjugacy for $\FFF$}
\label{sec:stab=qc}

The tangent family $\FFF$ is a holomorphic family  over $X=\C - \{0\}$. 
For this  family we prove: 

\begin{thm}
\label{thm:stab=qc} For the holomorphic family $\FFF$ over $\C - \{0\}$, the
J-stable parameters coincide with the quasiconformally stable parameters.
\end{thm}

\begin{proof}  Denote the asymptotic values of a  function in $\FFF$, by $v_1=
\lambda i$ and $v_2=-\lambda i$.  
If there exists an orbit relation for the singular values of a
function in $\FFF$, 
it has the form $$\fl{n}{\lambda}(v_i) =
\fl{m}{\lambda}(v_j) $$ 
 for some   integers $m,n>0$ and,  by
symmetry, for $i=1,2$.  

We claim that if there is an orbit relation for the singular values of a
function in $\FFF$, then the forward orbits of the asymptotic values are
finite.  If $i=j$ this is clear.  If $i\neq j$, 
 we may assume that  $n$ is even because if the relation holds for $m,n$, it
holds for $m+1,n+1$.  Suppose there is a relation, 
$\fl{n}{\lambda}(-\lambda i) = 
\fl{m}{\lambda}(\lambda i)$.  By proposition~\ref{thm:symmetry}, 
$\fl{n}{\lambda}(-\lambda i) = (-1)^{n}\fl{n}{\lambda}(\lambda i)
 = \fl{n}{\lambda}(\lambda i)$, so $\fl{n}{\lambda}(\lambda i) =
\fl{m}{\lambda}(\lambda i)$ and the forward orbits of both asymptotic 
values land on periodic cycles.

By the classification of stable behavior for functions in $\FFF$, the
Fatou set is non-empty if and only if there are either attracting
cycles, parabolic cycles or Siegel disks.   In all of these cases,
the forward orbits of the singular values are infinite.  Therefore,
if there are orbit relations, $J=\Chat$.  It follows that 
if $\lambda \in X^{stab}$ and
$J_{\lambda}\neq \Chat$ there are never any orbit relations so 
$\lambda \in X^{post}$.

On the other hand, if $\lambda \in X^{stab}$ and
$J_{\lambda}= \Chat$, there is a holomorphic motion defined on all
of $\Chat$ preserving the dynamics.  The  motion defines a topological
conjugacy that preserves the 
asymptotic values and their orbits so that 
$\lambda \in X^{post}$.

Thus,  whenever $\lambda \in X^{stab}$, 
$\lambda \in X^{post}$.    
\end{proof}

\subsection{Hyperbolic components and the density conjecture}
\label{sec:statement}

The hyperbolic maps form a natural subset of the J-stable maps.  
Set 
 $$\HHH = \{ \lambda \in \C-\{0\} : f_{\lambda}\;
\mbox{ has an attracting periodic cycle}\}.$$
Then the components of $\HHH$ are  components of  $X^{stab}$.

As an immediate corollary to theorem~\ref{thm:stab=qc} we have the
following result 
(which was proved independently in \cite{Jiang}):

\begin{cor}
\label{thm:theorem2} Each component  $\Omega$ of $\HHH$ is a component of
$X^{qc}$; 
  that is, for any ${\lambda}_1 $ and ${\lambda}_2 $ in
the same component of $\HHH$ there exists a quasiconformal map $\phi:
\hat\CC\to \hat\CC $ such that $ f_{{\lambda}_1}\circ \phi= \phi
\circ f_{{\lambda }_2}$.
 \end{cor}

It is natural to ask if 
 the hyperbolic components are the
only open components of the J-stable maps.  The  evidence of the
 computer pictures (figures~1,~2 and~3) indicate an affirmative answer.  

\begin{conj}[Density Conjecture]
\label{thm:density} For the holomorphic family $\FFF$ defined over $\C -
\{0\}$, 
the hyperbolic components are open and dense in the J-stable maps.
\end{conj}

\subsubsection{Invariant measurable line fields}

Let $Y$ be a component of the J-stable set $X^{stab}$ containing the
point $\lambda_0$ and let $\phi$ be a holomorphic motion of $\FFF$ over $Y$
with basepoint $\lambda_0$.  For each $\lambda \in Y$, the map
$\phi_{\lambda}:\Chat \rightarrow \Chat$ is quasiconformal.  Since $\phi$
respects the dynamics, the Beltrami
differential  $\mu_{\lambda}(z) = (\partial
\phi_{\lambda}/\partial{\bar{z}})/(\partial\phi_{\lambda}/\partial{z})$
satisfies
\begin{equation}
\label{eqn:beltinv} 
\mu_{\lambda}(f_{\lambda_0}(z))\frac{\overline{f_{\lambda_0}'(z)}}{f_{\lambda_0}'(z)}
= \mu_{\lambda}(z). 
\end{equation}
It  thus determines an {\em $f_{\lambda}$-invariant measurable line field on
$\C$}.

If $Y$ is hyperbolic, the support of $\mu_{\lambda}$ is contained in the
Fatou set of $f_{\lambda}$.  The measurable line field descends, via  the
grand orbit relation,  to the 
quotient of the Fatou set minus the closure of the grand orbits of the
singular values.  This quotient is 
either a pair of punctured tori of the same modulus or a single twice punctured
torus on which the punctures are symmetric.  (See \cite{McS} or
\cite{Jiang} for a fuller discussion). 

If $Y$ is not hyperbolic, $f_{\lambda}$ has no attracting cycles, and since
$\lambda$ is J-stable, $f_{\lambda}$ has no parabolic cycles or Siegel
disks.  Thus for any 
$\lambda$ in a non-hyperbolic component, $J_{\lambda}=\Chat$.

Now suppose that $J_{\lambda_0}=\Chat$ and $J_{\lambda_0}$ carries an
$f_{\lambda_0}$ measurable invariant line field $\mu$, $||\mu|| =1$.  For
any $t$, $|t|<1$, by the Measurable Riemann mapping theorem, \cite{Ahl-Bers},
$g_t(z) = \phi^{t\mu}\circ f_{\lambda_0}\circ(\phi^{t\mu})^{-1}$ is
meromorphic with exactly two asymptotic values and depends holomorphically
on $t$.  By
corollary~\ref{thm:nev}, normalized properly, $g_t=f_{\lambda}$ for some
$\lambda$ and $J_{\lambda}=\Chat$.  
Therefore, we obtain a  holomorphic motion of $\FFF$ defined over an open
disk in $X^{qc}$ based at ${\lambda_0}$. 

  By
theorem~\ref{thm:stab=qc}, 
conjecture~\ref{thm:density} is  thus equivalent to

\begin{conj}
\label{thm:density1} There is no disk $B$ in the set of quasiconformally
stable parameters for the family $\FFF$ such that the Julia set
$J_{\lambda}=\Chat$ for all $\lambda \in B$ and $J_{\lambda}$ carries an
invariant measurable line field.
\end{conj}

\section{The parameter space: combinatorial structure}

In this section  we will describe how the hyperbolic components fit together. 
  Our description
is analogous to the combinatorial description of the Mandelbrot set for
quadratic polynomials given by the periods of the attracting cycles.  For
example, as we saw in proposition ~\ref{thm:originattr} the component
$\Delta^* =\{| \lambda|<1\}$ has properties analogous to those of the
exterior of the Mandelbrot set.  The components of the Mandelbrot set have
a distinguished point, the center, at which the periodic cycle is
superattracting.  Because maps in $\FFF$ have no critical points, however,
the hyperbolic components do not have centers.  Nevertheless, we shall see
that they have a distinguished boundary point that we call a virtual
center. 

 The computer pictures (see figures 1, 2 and 3) of the
$\lambda$-plane drawn by W. H. Jiang suggest that all components of $\HHH$
except $\Delta^*$ appear in pairs and that each component pair has a unique
common boundary point.  We shall see that this point is the virtual center
of each component so we call it the {\em virtual center of the component
pair}.  The virtual center of the pair of unbounded components  is the
point at infinity.

In the next few sections we prove that the computer pictures of parameter
space in
figures 1, 2 and 3 are valid.

\begin{figure}
\label{fig:unitcirc}
\vspace{1.2in}
\centerline{\psfig{figure=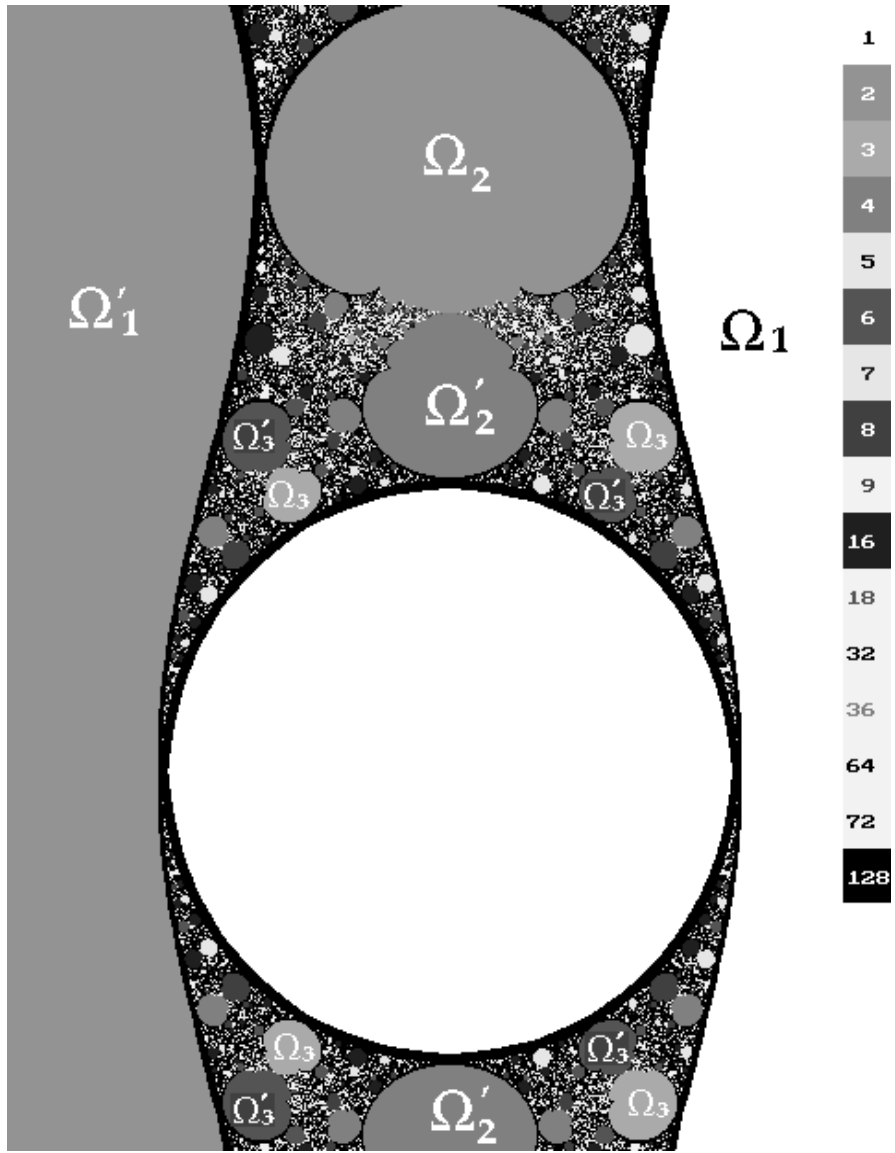,height=6in}}
\caption{The parameter plane near the unit circle.  The component pairs bud
off the unit circle at the endpoints of rational internal rays.}
\end{figure}

\begin{figure}
\label{fig:ball}
\vspace{1.2in}
\centerline{\psfig{figure=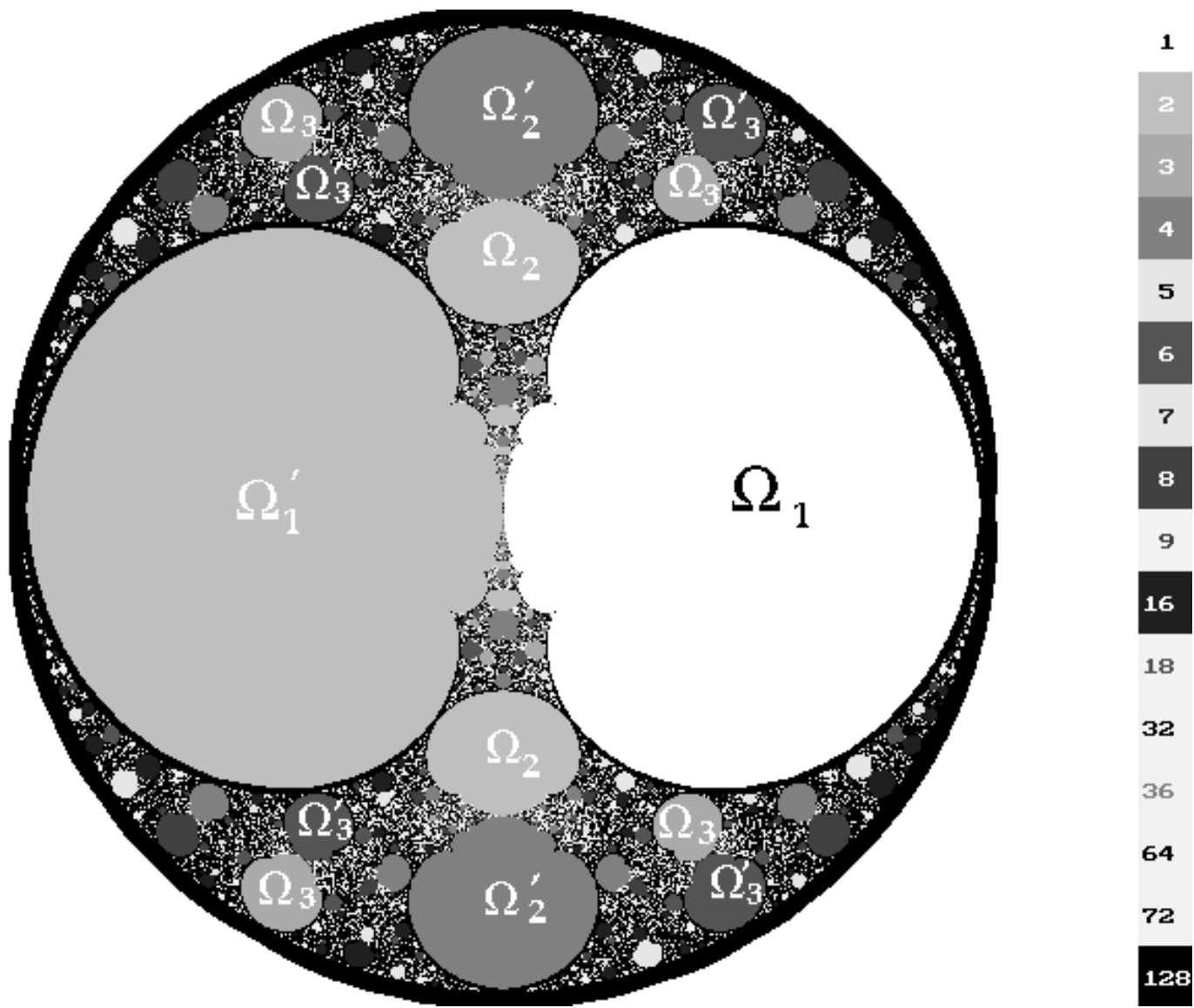,height=6in}}
\caption{The parameter plane inverted by $\lambda \mapsto 1/\lambda$.
The largest pair of components are the
$\Omega_{1},\Omega_{2}$ of corollary 5.4.}
\end{figure}

\begin{figure}
\label{fig:imagaxis} 
\vspace{1.2in}
\centerline{\psfig{figure=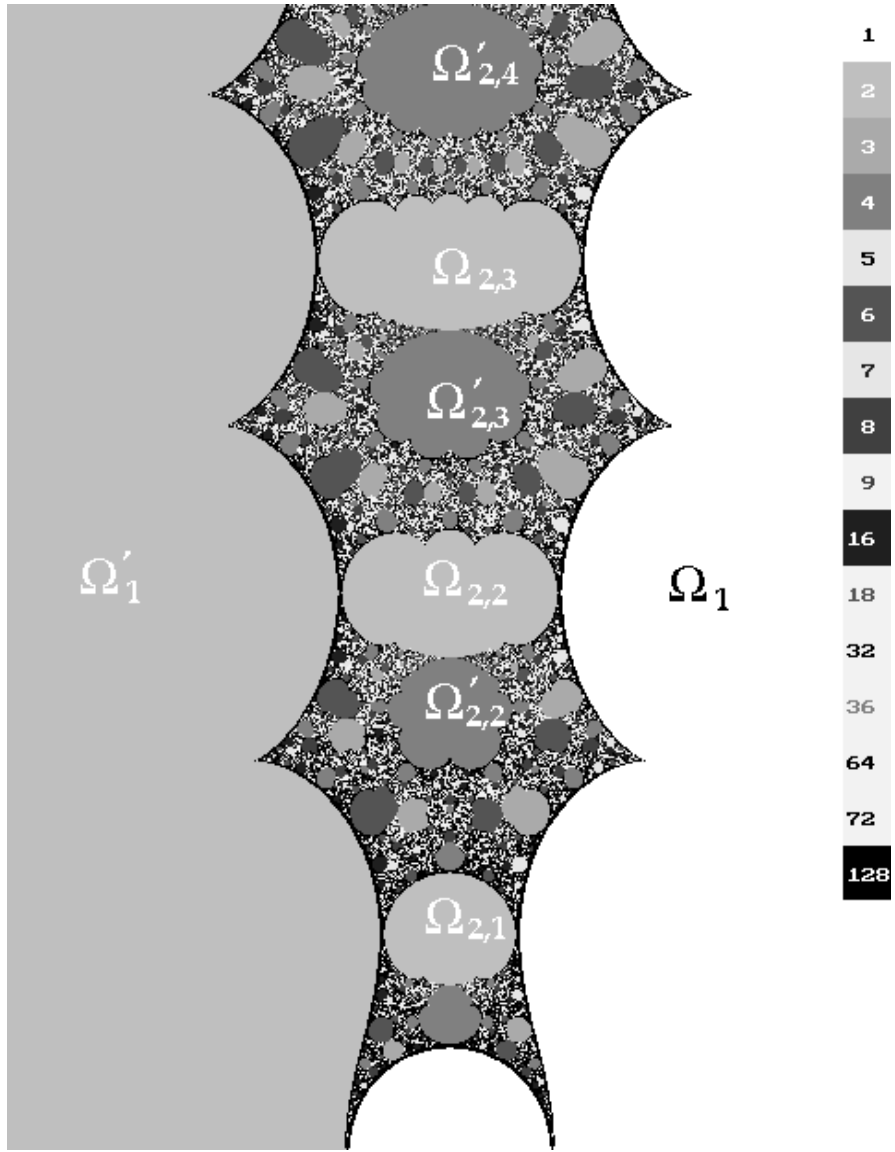,height=6in}}
\caption{The
$\lambda$ plane with the unit disk at the bottom of the figure. The
component pairs $\Omega_{2,n},\Omega_{2,n}'$ on the imaginary axis
with virtual center at $(n+1/2)\pi$ are clearly visible.}
\end{figure}

\subsection{The hyperbolic components}

By corollary~\ref{thm:theorem2}, 
it is  clear that the period of all maps in a given component is the
same. We therefore denote a generic hyperbolic component by $\Omega_p$
where $p$ is 
the period of the attracting cycle.

The following proposition is proved in \cite{DK1}.

\begin{prop}
\label{thm:reallambda}
Suppose $\lambda \in \RR$, $|\lambda| > 1$. Then
$J_{\lambda}=\RR $ and if
 $\lambda >1$, $f_{\lambda}$ has two attracting fixed points
$z_1=it$ and $z_2=-it$ for some $t=t(\lambda) \in \RR$, while if $\lambda
<-1$, $(z_1,z_2)$ is an attracting cycle of period $2$.
\end{prop}

Applying theorem~\ref{thm:qcdense} we obtain the corollary 

\begin{cor}
\label{thm:specialpair}
There is pair of components $ (\Omega_1,\Omega_2)$ containing the
positive and negative rays of the real axis with $|\lambda| >1$
respectively. For every $\lambda \in \Omega_i, \, i=1,2$,
$J_{\lambda}$ is a quasi-circle passing through $\infty$.
\end{cor}

This pair of  components can be characterized further by

\begin{cor}
\label{thm:bigcomp}
The component $\Omega_1$ containing the positive ray of the real axis,
$\lambda >1$, is the unique hyperbolic  component   such that
$f_{\lambda}$ has two distinct attracting fixed points; $\Omega_2$
containing the negative real ray is the unique component such that
$f_{\lambda}$ has a single period two cycle attracting both asymptotic
values. The boundary curve of $\Omega_1$ is asymptotic to the lines $t
\pm i e^{|2t|}$ as $t = \Re \lambda \to \infty$ and the boundary curve of
$\Omega_2$ is asymptotic to the lines $t \pm i e^{|2t|}$ as $t = \Re
\lambda \to -\infty$.
\end{cor}

\begin{proof} We need to show that if $|\lambda|>1$ and $f_{\lambda}$ has an
attracting fixed point, then $\lambda \in \Omega_1$.  Suppose $z =
\lambda\tan z$ and $|\lambda \sec^2 z| < 1$.  We can rewrite this condition
as $|h(u)| = |u/ \sin u| < 1$ where $u=2z$.  The locus $|h(u)|=1$ in the $u
= x + i y$-plane consists of four branches, meeting at the origin and
symmetric with respect to the origin.  The branches are asymptotic to the
lines $|x| \pm i e^{2|x|}$ as $|x| \to \infty$ and $|h(u)| <1$ in the upper
and lower half regions.  The function $ T(u)= u /(2 \tan 2u)$ satisfies
$|T(u)| \leq 1$ in a bounded set.  The complement of this set intersects
the upper and lower half regions in two connected, simply connected
unbounded domains.  The function $\lambda = T(u)$ is even and maps each of
these domains to a domain $\Omega$ in the $\lambda$ plane.  It is clear
that $\Omega$ must contain every $\lambda, |\lambda|>1$ for which
$f_{\lambda}$ has an attracting fixed point.  Moreover, $T$ maps the lines
$x = 0, y\neq 0$ to $\lambda >1$ so that $\Omega = \Omega _1$.

The asymptotic behavior of $\partial{\Omega_1}$ follows directly from that
of the curves $|h(u)|=1$.  By symmetry $\partial{\Omega_2}$ behaves the
same way.
   \end{proof}

In the next three propositions    we summarize  the results    on the
 properties of the remaining components of $\HHH$.  We omit
 the proofs of propositions~\ref{thm:prop2.13} and~\ref{thm:prop2.14}
 since they are straightforward.

\begin{prop}
\label{thm:prop2.13} Let $\Omega_p$ be a component of
$\HHH - \Delta^*$ such that for $\lambda \in \Omega_p$
$f_{\lambda}$ has an attracting cycle of period $p$. Then
either
\begin{itemize}
\item $f_{\lambda}$ has one attracting periodic cycle;
 both asymptotic values are in the immediate basin of the cycle
 and $p$ is even, or
\item $f_{\lambda}$ has two attracting periodic cycles of
 period $p$; both cycles have the same multiplier; the cycles
 are symmetric with respect to the origin and a basin of each
of them contains one of the asymptotic values. 
\end{itemize}
\end{prop}

\noindent{\bf Definitions:} As usual, 
let $m_{\lambda}$ denote the multiplier of an attracting or neutral
periodic cycle of $f_{\lambda}$ containing the point $z_{\lambda}$.  If
$\Omega_p$ is an arbitrary hyperbolic component of $\HHH$ and $\Delta^*$
is the unit disk punctured at the origin, the {\em eigenvalue map}
$m:\Omega_p \rightarrow \Delta^*$is defined by $\lambda \mapsto
m_{\lambda}$. For each $\alpha \in \RR$ the {\em internal ray}
$R(\alpha)$ is defined by $R(\alpha)=m^{-1}(r e^{2\pi i \alpha}), 0 < r < 1$.

\begin{prop}
\label{thm:prop2.14}
 The eigenvalue map $m: \Delta^* \rightarrow \Delta^*$ is the
identity.
 For each component $\Omega_p$ of $\HHH - \Delta^*$ the
eigenvalue map is an infinite degree
 regular covering map.
\end{prop}

The eigenvalue map lifts to a conformal
isomorphism
$\tilde{m}$ of $\Omega_p$ onto the upper half plane $\HH$ defined by
$$\tilde{m}(\lambda)\mapsto 2 \pi
\arg{m(\lambda)}-i\log|m(\lambda)|$$
where the branch of the logarithm is chosen so that the internal ray $
R(0)$ is mapped to
$-i \log r$.
 Under this map the boundary of $\Omega_p$
corresponds to the real axis $\RR$ and the common endpoint of
all the internal rays of $\Omega_p$ is a boundary point
 that corresponds to the point at infinity under $\tilde{m}$.

\Definition The boundary point corresponding to the point
 at infinity of the component $\Omega_p$ is called its
{\em virtual center}. Recall that the hyperbolic components have
 no center since the functions in $\FFF$ have no critical points.
The virtual center for the special components  $\Omega_1$ and
$\Omega_2$ above is $\infty$.

\subsection{Virtual centers and component pairs}

From now on we reserve the notation $\Omega_p$ for the components of
$\HHH - \Delta^*$ with two distinct attracting cycles.  We use
$\Omega_p'$ for the components with a single cycle of period $2p$.
In particular, the component called $\Omega_2$ above  is now renamed
 $\Omega_1'$.

We will show that if $p>1$, the asymptotic values of the functions
corresponding to the virtual centers of either $\Omega_p$ or
$\Omega_p'$ are prepoles of order $p-1$. For $p=1$, the common virtual
center of $\Omega_1$ and $\Omega_1'$ is infinity.  The proof has two
parts: we prove the statement first under the assumption that the
components are bounded. We remove this assumption in the next section
by describing the deployment of the $\Omega_2$ components.

\begin{prop}
\label{thm:jiang} For any bounded hyperbolic component
 $\Omega_p$ or $\Omega_p'$  with $p>1$, the virtual center $\lambda^*$ is
finite 
and $\fl{p-1}{\lambda^*}(\lambda^* i)= \infty$; that is $\lambda^* i$
is a prepole of order $p-1$.
\end{prop}

\begin{proof}
 For $\lambda \in \Omega_p$, let $z_0=z_0(\lambda)$ be the attracting
periodic point of $f_{\lambda}$ of period $p$ such that $z_0$ belongs to
the component $D_0 $ of the regular set that contains the asymptotic tract
and such that $\lambda i$ and $z_1$ are both in $D_1=f(D_0)$.  Denote the
preimage of $z_0$ in the periodic cycle by $z_{p-1}$ and the preimage of
$D_0$ by $D_{p-1}$; then, for some $n$, $z_{p-1}=f^{-1}_n (z_0)$.  Since
$p>1$, the domains $D_{p-1}$ and $D_0$ are different and the map $f:D_{p-1}
\to D_0$ is bijective.  Since $D_0$ contains the asymptotic tract and the
map is bijective, there must be a unique pole, $s_n$, on
$\partial{D_{p-1}}$.  To see this, note that there is a preasymptotic tract
at $s_n$ in $D_{p-1}$ containing a preimage of either $z=iy$ or $z=-iy$ for
large $y>0$.  So if $\partial{D_{p-1}}$ contained any other pole, there
would be a preasymptotic tract in $D_{p-1}$ at this pole containing a
preimage of the same segment and $f|_{D_{p-1}}$ would not be injective.

Since $\Omega_p$ is a component of the J-stable set, 
the Julia set and hence the prepoles vary continuously as we vary
$\lambda$. 
Since the poles are fixed however, $s_n$ remains the unique pole on
$\partial{D_{p-1}}$.

Suppose that  $\lambda$ moves along the internal ray $R(\alpha)$ in
$\Omega_p$
to the virtual center $\lambda^*$ as $r \rightarrow 0$, so that
 $$\lim_{\lambda
\stackrel{R}\rightarrow \lambda^*}
 m_{\lambda}=0.$$
  Since $(\tan z)^{\prime}=\sec^2 z$, and $\lambda=z_{i}/\tan z_{i-1}$,
the
multiplier $m_{\lambda}$ can be written as
$$m_{\lambda}= [\fl{p}{\lambda}(z_0({\lambda}))]^{\prime}=
\Pi_{k=1}^p f^{\prime}_{\lambda}[\fl{k-1}{\lambda}(z_0({\lambda}))]=
 \Pi_{k=1}^p 2z_i/\sin 2z_{i-1}.$$

The only way some factor may tend to $0$ as $\lambda \to \lambda^*$ is for
$\sin 2z_{i-1} \to \infty$ for some $i$, or equivalently, for $\Im z_{i-1}
\to \infty$.  Since $z_0$ is in the asymptotic tract, we conclude $\Im z_0
\to \infty$.

 By  hypothesis  $p>1$ and  $\lambda^*\neq \infty$ so
$z_{p-2} \neq z_{p-1}$.
  We have
$$\lim_{\lambda\stackrel{R}\rightarrow\lambda^*}\lambda \tan
z_{p-1}(\lambda)=
\lim_{\lambda \stackrel{R}\rightarrow \lambda^*} z_0 = \infty.$$

We conclude  further that
$\lim_{\lambda\stackrel{R}\rightarrow
\lambda^*}z_{p-1}(\lambda)=\lim_{\lambda\stackrel{R}\rightarrow
\lambda^*}f^{-1}_{\lambda,n}(z_0(\lambda))=s_n$ and
$\lim_{\lambda\stackrel{R}\rightarrow \lambda^*} z_1(\lambda)=\lambda^* i$
so that $\lambda^* i$ is a prepole of
of $f_{\lambda^*}$ of order $p-1$. The other periodic orbit behaves
symmetrically and hence $-\lambda^*i$ is also a prepole.

If  $\lambda \in
\Omega_p'$, $f_{\lambda}$ has a single cycle containing the points $z_0$ and
$z_p=-z_0$.  These lie in symmetric components $D_0$ and $D_p$
containing the asymptotic tracts.  We argue as above that there are
unique poles $s_n$ and $-s_n$ on their respective boundaries. If $\lambda^*$
is the virtual center, we again
conclude that
$\lim_{\lambda\stackrel{R}\rightarrow \lambda^*} z_1(\lambda)=\lambda^* i$
and also that
$\lim_{\lambda\stackrel{R}\rightarrow \lambda^*}
 z_{p+1}(\lambda)=-\lambda^* i$  so that
 $\lambda^* i$ is a prepole of
of $f_{\lambda^*}$ of order $p-1$.
  \end{proof}

\medskip 

The techniques developed by Douady, Hubbard and Sullivan for quadratic
polynomials, showing that the boundary of a hyperbolic component is a
piecewise analytic curve adapt easily to show that $\partial\Omega_p$
and $\partial\Omega_p'$ are also piecewise analytic curves
(\cite{Carlesson,DH,Jiang}).
The main difference is that the eigenvalue map is an infinite degree
universal cover.
  Similarly, it is straightforward to
modify the local techniques used to describe the root and bud
structure of the Mandelbrot set (\cite{Jiang}).
Note the following however: for $\lambda \in \Omega_p'$, both $z_0$
and $-z_0$ belong to the attracting periodic cycle, $-z_0=f^p(z_0)$
and $m_{\lambda}=\lambda^{2p}(\Pi_{i=0}^{p-1}\sec^z_i)^2$.  We thus
redefine the map from the universal cover of $\Omega_p'$ to the upper half
plane $\HH$ by $\tilde{m}\to 1/2(\arg m_{\lambda} -i \log|m_{\lambda}|)$.
 We omit the proofs of the following propositions (see \cite{Jiang}).

\begin{prop}
\label{thm:budcomp} Let $q>1$ and let $\lambda \in \partial\Omega_p$ be such
 that $m(\lambda)$ is a primitive $q$-th root of unity. Then there is a
``bud'' component $\Omega_{pq}$ with two attractive cycles of order $pq$
tangent to $\Omega_p$ at $\lambda$. Similarly, if $\lambda \in
\partial\Omega_p'$, and $m(\lambda)$ is a primitive $q$-th root of unity,
there is either a bud component $\Omega_{pq}$ or a bud component
$\Omega_{pq}'$ attached at $\lambda$.  \end{prop}

\begin{prop} At points of $\partial\Omega_p$ and $\partial\Omega_p'$
such that $\tilde{m}(\lambda)= 2 \pi n$, $n \neq 0$, there is either a cusp
or there may be  a ``root'' component $\Omega_d$
where $d/p$.  If $p=1$, the root component is $\Delta^*$.  There are at
most finitely many roots.
\end{prop}

 At a cusp of $\partial\Omega_p$, there are two cycles; at each point of
each periodic cycle there is a single petal that contains the forward orbit
of the asymptotic value attracted to that cycle.  At a cusp of
$\partial\Omega_p'$ there is a single cycle; at each periodic point again
there is a single petal, this time containing the forward orbits of both
asymptotic values.

At a point of either $\Omega_p$ or $\Omega_p'$ where a bud bifurcation
of order $q>1$ occurs the points of the parabolic cycle always have
$q$ petals. If the point is on $\Omega_p$ there are two cycles; at each
parabolic periodic point there are $q$ petals and the forward orbit
of the asymptotic value attracted to the cycle cycles through all
the petals.
 If the bud point is
on $\Omega_p'$, there is a single cycle of period $2p$ and again at each
point in the cycle there are $q$ petals.  If the
forward orbits of  asymptotic values belong to distinct sets of petals,
the bud component has two distinct attractive cycles of order $pq$,
 while if the forward orbits of asymptotic values belong to the same
set of petals the bud component has a single cycle of period $2pq$.
If the bud point is on $\Omega_p'$, each petal contains the orbit of a
different asymptotic value.  On $\Omega_p'$ there are bifurcation
points with $q=1$; this happens if at each parabolic point there are
two petals and the forward orbits of the asymptotic values are in
different cycles of petals.

\subsection{Components Tangent to $\Omega_1$ and $\Omega_1'$}

It follows from proposition~\ref{thm:budcomp} that there is a sequence
$\lambda_k \in \partial\Omega_1$ with
$\tilde{m}(\lambda_k)=(2k+1)\pi$, $k \in \ZZ$ and a sequence of bud
components $\Omega_{2,k}$ tangent to $\partial\Omega_1$ at
$\lambda_k$.  We know from theorem~\ref{thm:dichot} and
 proposition~\ref{thm:reallambda} that $\Omega_{2,k} \cap \RR = \emptyset$.
Thus below we will fix $k$, set $\Omega_2 =\Omega_{2,k}$
 and  assume  that $\Im \lambda > 0$
for  $\lambda \in \Omega_{2,k}$.

We will prove that the virtual center of the component
$\Omega_{2}=\Omega_{2,k}$ is $s_k i$ and so is finite;
 we will use this fact to
conclude that, for all $p>1$, the components $\Omega_p$ and $\Omega_p'$
are bounded.  We need
the following lemmas.

\begin{lemma}
\label{thm:boundedvirt1}  Let $\Omega_{2}$ be a  hyperbolic
 component such that for $\lambda\in \Omega_{2}$, $f_{\lambda} $ has
two period two attracting cycles $z_0,z_1$ and $-z_0,-z_1$. If the
virtual center $\lambda^*=\infty $, then   
for $j=0,1$, as $\lambda$ varies along some internal ray $R(\alpha)$ in
$\Omega_{2}$, 
$$\pm z^*_j=
\lim_{\lambda\stackrel{R}\rightarrow \lambda^*}
(\pm z_j({\lambda}))=\infty. $$
\end{lemma}

\begin{proof}  For $\lambda \in \Omega_{2}$  let
 $z_0=z_0(\lambda)$ belong to the component $D_0 $ of the regular set
that contains the asymptotic tract of $\lambda i$; then $\lambda i$
and $z_1$ are both in $D_1=f_{\lambda}(D_0)$.  The other periodic
orbit $-z_0,-z_1 $ behaves symmetrically.  Let us suppose that $\lambda$
moves along the internal ray $R(\alpha)$ in $\Omega_{2}$ to a limit
point $\lambda^*=\infty$.  Then $\lim_{\lambda\stackrel{R}\rightarrow
\lambda^*}m_{\lambda}=0$.  It follows (compare the proof of
proposition~\ref{thm:jiang}) that $\Im z_0(\lambda)\to +\infty $.
  We want to prove that $
z^*_1=\lim_{\lambda\stackrel{R}\rightarrow \lambda^*}z_1({\lambda})
=\infty$. If not, there exists a sequence $\lambda_n
\stackrel{R}\rightarrow \lambda^*=\infty$, such that
$\lim_{\lambda\stackrel{R}\rightarrow \lambda^*}z_1({\lambda_n}) = c
\neq\infty$.
Since  $z_1(\lambda_n)=\lambda_n\tan z_0 (\lambda_n)$,
 we must have
$\lim_{\lambda_n\stackrel{R}\rightarrow \lambda^*}\tan(z_0(\lambda_n))
=0$. Thus either the curve $z_0(\lambda)$ is bounded and $
\lim_{\lambda\stackrel{R}\rightarrow \lambda^*}(z_0({\lambda}))=m\pi $
for some integer $ m$, or the curve $z_0(\lambda)$ is unbounded but
comes arbitrarily close to infinitely many integral multiples of $\pi$.
Either possibility contradicts $\Im z_0(\lambda)\to
\infty$.  
\end{proof}

\begin{lemma}
\label{thm:per2bdd}
Let $\Omega_{2}$ be the bud component tangent to $\Omega_1$ at
$\lambda_k$ as above.  The virtual center $\lambda^*$ is
equal to $s_k i$.
\end{lemma}

\begin{proof} We claim that $\lambda^*$ is finite.  If not,
by lemma~\ref{thm:boundedvirt1},  taking the limit
along an internal ray, the periodic points $z_0,z_1$
 go to infinity,
 and in particular, $\Im z_0$ goes to infinity.
We claim this cannot happen.  For readability we suppress the
dependence on $\lambda$ and on the ray; we
 set $\lambda = \lambda_1 + i \lambda_2$,
 and $z_j= x_j + i y_j$, $j=0,1$.  Since
$z_1=\lambda\tan z_0$,
\begin{equation}
\label{eqn:x0}
x_1=\frac{ \lambda_1\sin(2x_{0})-
\lambda_2\sinh(2y_{0})}{\cos(2x_{0}) +
\cosh  (2 y_0) }.
\end{equation}

and
\begin{equation}
\label{eqn:y0}
y_1=\frac{ \lambda_1\sinh(2y_{0})+
\lambda_2\sin(2x_{0})}{\cos(2x_{0}) +
\cosh  (2 y_0) }.
\end{equation}

From lemma~\ref{thm:boundedvirt1} we have $y_0 \to \infty$ and by
corollary~\ref{thm:bigcomp} we have $\lambda_2 \geq e^{|2\lambda_1|}$ and
$\lambda_2 \to \infty$.  From equation~\ref{eqn:x0} we see that $|x_1|
\approx \pm \lambda_2 \to \infty$.  From the periodicity we obtain formulas
for $x_0$ and $y_0$ by interchanging $0$ and $1$ in equations~\ref{eqn:x0}
and~\ref{eqn:y0}.  Since $|x_1| \to
\infty$ the term $\lambda_2 \sin(2x_1)$ in the equation for $y_0$
oscillates so $|\lambda_1 \sinh(2y_1)|$ must grow faster than  $|\lambda_2
\sin(2x_1)|$; this implies $2y_0 \approx \lambda_1 \to \pm\infty$.  The
periodicity
again implies $|x_0| \approx \pm \lambda_2 \to \infty$ and $2y_1 \approx
\lambda_1 \to \pm\infty$.  To estimate the
multiplier of the cycle we have $$ |\lambda \sec^2 (z_i)| \approx |\lambda_2|
e^{\pm \lambda_1},$$ so that either $$|m(\lambda)|= \lambda_2^2 e^{\pm 2\lambda_1}
\quad {\mbox or }\quad |m(\lambda)|=\lambda_2^2.$$
Since the cycles of $f_{\lambda}$ are $(z_0,z_1)$ and $(-z_0,-z_1)$, the
cycles of $f_{-\lambda}$ are $(z_0,-z_1)$ and $(z_0,-z_1)$ and we conclude
that either $|m(\lambda)|$ or $|m(-\lambda)|$ grows with $\lambda_2^2$ and so
cannot tend to zero.

Thus, $\lambda^*$ is finite and the asymptotic value of
$f_{\lambda^*}$ must be a prepole of order 1 and so
must be $s_n i$ for some pole $s_n$. We may draw a curve
$\gamma_k$ in $\Omega_{2,k}$ from the bifurcation point
$\lambda_k$ to the virtual center $s_n i$. These curves
are disjoint and hence occur in order. We shall see in
theorem~\ref{thm:centers} that each $s_n i$ is a virtual center
of some component so the endpoint of $\gamma_k$ must be
$s_k i$.
 \end{proof}

If $\lambda$ is the bifurcation point on $\partial\Omega_1$ for the
bud $\Omega_2$, the point $-\bar{\lambda}$ is a bifurcation point on
$\partial\Omega_1'$ for a bud component $-\bar\Omega_2$.  We leave it
to the reader to check that the single period two cycle bifurcates to
two period two cycles.  The point here is that the forward orbits of
the asymptotic values are in distinct cycles of petals at both
bifurcation points.

\begin{prop} The hyperbolic components $\Omega_p$ and $\Omega_p'$
are bounded for all $p>1$.
\end{prop}

\begin{proof}  For each integer $k$ we can draw the curves $\gamma_k$
and $-\bar{\gamma}_k$ in the components $\Omega_{2,k}$ and
$-\bar{\Omega}_{2,k}$.  Any hyperbolic component, not equal
to $\Omega_1, \Omega_1',\Omega_{2,k}$ or $\bar\Omega_{2,k}$
must be contained inside a region bounded by curves $\gamma_k,
-\bar\gamma_k, \bar\gamma_k, -\gamma_k$ and arcs of
$\partial\Omega_1$ and $\partial\Omega_1'$.
 \end{proof}

\subsection{Combinatorial invariants}

If $\lambda^*$ is a virtual center, we know that
$J_{\lambda^*}=\hat\CC$, but the two sets
$$\{\infty, \pm \lambda^* i,
f_{\lambda^*}(\pm \lambda^*i),
\ldots,\fl{p-2}{\lambda^*}(\pm \lambda^* i)\}$$ are cycles,  considered as
appropriate limits.   They have a property
that superattractive cycles have, namely, that each contains a singular
value, the asymptotic value. 

Set
 $$\CCC_p = \lbrace \lambda: f^{p}_{\lambda}(\lambda i)
 = \infty \rbrace, \,\, \CCC=\cup_{1}^{\infty}\CCC_p.$$

The next
theorem implies that if $\lambda \in \CCC_{p-1}$ there
is a pair
 of components $(\Omega_p,\Omega_{p}')$ with $\lambda$ as virtual
center.  It follows that the itineraries of the virtual centers give a
combinatorial description of the component pairs.

\begin{thm}
\label{thm:centers}
Let ${\lambda}_0i$ be a prepole of order $p-1$ with
$\f{p-2}({\lambda}_0 i)=s_n $. Then ${\lambda}_0$ is the virtual
center of a  component pair $ (\Omega_p,{\Omega}_{p}')$ and
$\lambda_0 \in\partial{\Omega}_p \cap \partial{\Omega}_{p}'$,
where for $\lambda \in \Omega_p, f_{\lambda} $  has two attracting
cycles of period $p$ and for $\lambda \in \Omega_{p}',f_{\lambda}$
  has one  attracting cycle of period $2p$.
  \end{thm}

\begin{proof}
To prove this theorem, we need to show that there exist $\lambda$ and
$\lambda'$ arbitrarily close to $\lambda_0$ such that $f_{\lambda}$ has an
attracting periodic cycle of period $p$ and $f_{\lambda'}$ has an attracting
periodic cycle of period $2p$.  To construct the attracting periodic cycle
we want to find $\lambda$ and a domain $\TT$ in the dynamic plane such that
$\f{p}_{\lambda}(\TT) \subset \TT$.

For any choice of $r>0$ we can choose asymptotic tracts
 $\AAA^+=\{\Im z >r\}$ and $\AAA^-=\{\Im z < -r\}$ for
$\pm \lambda_0 i$ respectively. Since $\f{p-2}(\lambda_0 i)=s_n$,
we have $\f{p-2}(-\lambda_0 i)=-s_n$.  Set
 $U = f_{\lambda_0}(\AAA^+)$ so that $U$ is a neighborhood
of $\lambda_0 i$ and let $V$ be the corresponding neighborhood
 in the parameter plane; that is, $\lambda \in V$ if and only if
$\lambda i \in U$.

For $\lambda \in V$ take the common preasymptotic
 tracts
$$\BBB_n^{\pm}= \cap_{\lambda \in V} f^{-1}_{n,\lambda}(\AAA^{\pm})$$
attached to the pole $s_n$. We can find $0 <\delta = \delta(r)$
 such that  $|\arg(\lambda -
\lambda_0)| < \delta $ if
 $\lambda \in V$. Hence the
angle between $f^{-1}_{n,\lambda}(\RR)$ and $\RR$ at $s_n$ is
bounded and $\BBB^{\pm}$ contains some triangular domain  with
 one vertex at $s_n$. Define a map $g:V \rightarrow \CC$ by
$g(\lambda)= \fl{p-2}{\lambda }(\lambda i)$ for $\lambda \in V$.
Then $g(V)$ is an open set containing the pole $s_n$ and there
 exist open sets $V^{\pm} \subset V$ such that
 $V^{\pm} = g^{-1}(\BBB_n^{\pm})$. For any $\lambda \in V^+$,
$\fl{p-1}{\lambda}(\lambda i)$ belongs to an asymptotic tract
$\AAA^+=\{ \Im z > r'\}$ where possibly $r'<r$. Moreover, for
 inverse branches such that
 $$\fl{-(p-2)}{{\bf n}_{p-2},\lambda_0}(s_n)=\lambda_0 i,$$
 we have the property that
 $v_{\lambda}=\fl{-(p-2)}{{\bf n}_{p-2},\lambda}(s_n)\neq \lambda i$,
 so that there are  preimages
 $w_{k,\lambda} = f^{-1}_{k,\lambda}(v_{\lambda})$ in the
 upper half of each strip $L_k, k \in \bf Z$. These $w_{k,\lambda}$ depend
 continuously on $\lambda$ and if
$\lambda \rightarrow \lambda_0$, then
$\Im w_{k,\lambda} \rightarrow \infty$.

Let
$\eta_{\lambda}= |\lambda i - v_{\lambda}|$ and consider the
 ball $B_{\lambda}=B(v_{\lambda},\eta_{\lambda})$. Then
 $\fl{p-2}{\lambda}(B_{\lambda})$ is an open set containing $s_n$.
 Using the principal part of $f_{\lambda}$ we see that
 $\fl{p-1}{\lambda}(B_{\lambda})=\hat\CC - B(0,R_{\lambda})$
where $R_{\lambda} \approx |\fl{p-1}{\lambda}(\lambda i)|$
and $R_{\lambda} \to \infty$ as $\lambda \to \lambda_0$.
We claim that $\Im \fl{p-1}{\lambda}(\lambda i) > \Im w_{k,\lambda}$.
 Let
$$K=\max_{z \in \bar{U},\lambda \in
\bar{V}}|(\fl{p-2}{\lambda})^{\prime}(z)|.$$
Since $|\lambda|>1$ and $|\lambda | \gg \eta_{\lambda}$, we argue as in
the proof of proposition ~\ref{thm:proppp} and obtain
$$ \Im w_{k,\lambda}=\Im f_{\lambda}^{-1}(v_{\lambda})=\Im
f_{\lambda}^{-1}(\lambda i + \eta_{\lambda})
= \Im\left[ \frac{1}{ 2i} \log \frac{i\eta_{\lambda}}{2 \lambda
- i\eta_{\lambda}}\right].$$
Thus,
 $\Im w_{k,\lambda} \approx \frac{1}{2}|\log |\eta_{\lambda}||$.

Now set
$y_0=\frac{1}{2} |\log \eta_{\lambda}|$. Then
$|f^{-1}_{n,\lambda}(i y_0)-s_n|\approx 2\lambda/|\log \eta_{\lambda}|$
 and consequently
$$
|\fl{-(p-1)}{{\bf n}_{p-1},\lambda}(i y_0) - \lambda i|
 \geq |f^{-1}_{n,\lambda}(i y_0) -s_n|\cdot
\min_{w \in \fl{p-2}{\lambda}(U)}|
(\fl{-(p-2)}{{\bf n}_{p-2},\lambda})^{\prime}(w)|$$
and approximating,
$$|\fl{-(p-1)}{{\bf n}_{p-1},\lambda}(i y_0) - \lambda i| \geq {2
\lambda\over {K |\log \eta_{\lambda}|}}.$$
But $ |\lambda i- v_{\lambda}| = \eta_{\lambda}$ and since
 $\eta_{\lambda}$ is assumed small
$$\frac{2 \lambda}{K |\log \eta_{\lambda}|} \geq \eta_{\lambda}.$$
 Therefore
$$|\fl{p-2}{\lambda}(\lambda i)-s_n|\leq |f^{-1}_{n,\lambda}(i y_0)-s_n|$$
so that $$\Im \f{p-1}_{\lambda}(\lambda i) > i y_0$$
and $$\Im \fl{p-1}{\lambda}(\lambda i) > \Im w_{k,\lambda} \approx  y_0$$
as claimed.

Now we are ready to construct the domain $\TT$ inside an asymptotic
tract $\AAA^+$ for some $\lambda$ in $V$ with
$\fl{p}{\lambda}(\TT) \subset \TT$  proving that
$\lambda \in\Omega_{p}$. Set
 $r= r_{\lambda}= (1/2)|\log \eta_{\lambda}|-\epsilon$, and choose an
 asymptotic tract
 $\AAA^+=\lbrace z:\Im z > r_{\lambda}\rbrace$ so that
 $v_{\lambda} \in f_{\lambda}(\AAA^+)$. Let $I^{\pm}$ be two rays
 meeting at $s_n$  such that the triangular domain $T$ between
them containing the vertical direction is contained in
$f^{-1}_{ n,\lambda}(\AAA^+)$ and such that
$\fl{p-2}{\lambda}(\lambda i) \in T$. Let $\TTT$ be the
triangular region with vertex at $v_{\lambda}$ bounded by
${\cal I}^{\pm} = f^{-(p-2)}_{{\bf n}_{p-2},\lambda}(I^{\pm})$
and an arc of the boundary of $f_{\lambda}(\AAA^+)$ so that
 $\lambda i \in \TTT$. Finally, set
$\TT =\cup_{k \in \ZZ}f^{-1}_{k,\lambda}(\TTT)$. Then $\TT$
is an asymptotic tract whose ``horizontal'' boundary is
scalloped by preimages of $\partial \TTT$; that is, it is made up
 of arcs  $f^{-1}_{k,\lambda}({\cal I}^{\pm})$ that meet at
 $w_{k,\lambda}$ and are joined by intervals in the line
 $y=r_{\lambda}$.

Now consider
 $\tilde\TT=\fl{p-1}{\lambda}(\TTT)$;
this is a triangle with a vertex at infinity; the sides meeting there are
 rays and the third side is an arc of a circle centered at the origin with
 radius
slightly smaller than $|\fl{p-1}{\lambda}(\lambda i)|$.
 Because
$|\fl{p-1}{\lambda}(\lambda i)| > (1/2)|\log \eta_{\lambda}|$,
 we may assume
$|\partial{\fl{p}{\lambda}(\AAA^+)}| >(1/2)|\log \eta_{\lambda}|$
  for small $\epsilon$.  To claim that
$\fl{p}{\lambda}(\TT) \subset \TT$ we need to check two
conditions:
\begin{enumerate}
\item $\Im \fl{p-1}{\lambda}(\lambda i) > r_{\lambda}$ and
\item $f_{\lambda}(I^{\pm}) \subset \AAA^+$.
\end{enumerate}
 Now $\lambda$ was chosen so that
$\fl{p-1}{\lambda}(\lambda i)$ belongs to its asymptotic tract
hence changing the argument of $\lambda$ if necessary we can
insure  that 1. holds. We can insure 2. by decreasing the angle
between $I^{\pm}$ if necessary.

 In the above we chose
 $\lambda \in V^{+}$ such that
$\fl{p-2}{\lambda}(\lambda i) \in \BBB^+$. Next we show that
if $\lambda \in V^-$ and $\fl{p-2}{\lambda}(\lambda i) \in \BBB^-$
then $\lambda \in \Omega_{2p}$. Take inverse branches such that
$$\fl{-(p-2)}{{\bf n}_{p-2},\lambda_0}(s_n )=
\lambda_0 i \mbox{ and } \fl{-(p-2)}{{\bf n}_{p-2},\lambda_0}(-s_n )
= -\lambda_0 i.$$
Since $\fl{p-2}{\lambda}(\lambda i) \in \BBB^-$ we have
 $\fl{p-1}{\lambda}(\lambda i) \in \AAA^-.$ Now
$\fl{p}{\lambda}(-\lambda
i) \in \BBB^+$ so
$\fl{2p-1}{\lambda}(\lambda i) \in \AAA^+$. Therefore,  arguing
analogously to the above we can prove that $\lambda \in\Omega_{p}'$.
 \end{proof}

The following corollaries are immediate.
\begin{cor}
Assume  $p \geq 2$ and $\lambda \in \Omega_p \cup \Omega_p'$.
\begin{itemize}
\item if the asymptotic value $\lambda i$ is contained in the
preasymptotic tract $\BBB_v^+$ of the prepole $v$ of order $p-1$
(respectively $-\lambda i \in \BBB_v^-$),  then
 $f_{\lambda}$ has two attracting periodic cycles of period $p$ and
 $\lambda \in \Omega_p$.
\item if the asymptotic value $\lambda i$ is contained in the
preasymptotic tract $\BBB_v^-$, then there is a single periodic
 cycle of period $2p$ attracting both singular values and $\lambda \in
 \Omega_p'$ .
\end{itemize}
\end{cor}

\begin{cor} Let $(\Omega_p,\Omega_{p}')$ be a pair of hyperbolic
 components meeting at the virtual center $\lambda^*$. Suppose
$\lambda(t),\, t \in [0,1]$ is a curve in the parameter plane
 passing from $\Omega_p$ to $\Omega_{p}'$ through $\lambda^*$.
Then there is a {\em period doubling bifurcation} of attracting
periodic points of $f_{\lambda}$.
\end{cor}

\subsection{Duality of  the virtual centers and prepoles}

In this section we show that the virtual centers in the parameter
 plane play a role that is dual to the prepoles in the dynamic
plane. 
\begin{prop}
\label{thm:centersacc} Let $\lambda_n \in \CCC_{p-1}$ where
$\fl{p-2}{\lambda_n}(\lambda_n i) = s_n$; that is, $\lambda_n$
is the virtual center of a component pair $(\Omega_p,\Omega_{p}')$.
Then there is a sequence of component pairs
 $(\Omega_{p+1,k},\Omega_{p+1,k}'),\; k \in \ZZ$, with centers
 $\lambda_k \in \CCC_p $ where
 $\fl{p-1}{\lambda_k}(\lambda_k i)=
 s_k$ and $\lambda_k \rightarrow \lambda_n$ as
 $|s_k| \rightarrow \infty$.
\end{prop}

\begin{proof} Let $B(\lambda_n i,\epsilon)$ be an $\epsilon$
neighborhood  of the asymptotic value
$\lambda i$ of $f_{\lambda}$ in the dynamic plane. Let $D(\lambda_n,
\epsilon)$ be the corresponding
neighborhood in the parameter plane.  By assumption
$\fl{p-2}{\lambda_n}(\lambda_n i) = s_n$. Define a map
$g:D(\lambda_n,\epsilon) \rightarrow \hat\CC$ by
$g(\lambda)=\fl{p-1}{\lambda}(\lambda i)$. Then
$g(D(\lambda_n,\epsilon))$ is  a neighborhood of infinity.
 Therefore, for large $k$, there is a parameter $\lambda_k \in
 D(\lambda_n,\epsilon)$ such that
 $$ g({\lambda}_k) =  \fl{p-1}{\lambda_k}(\lambda_k i) =s_k. $$
 Thus $\lambda_k \in \CCC_p$ and by theorem ~\ref{thm:centers} $\lambda_k$
is  the   virtual   center  of   a   component pair
$(\Omega_{p+1,k},\Omega_{p+1,k}')$; moreover,
$\lambda_k \rightarrow \lambda_n$ as
$|s_k| \rightarrow \infty$.
 \end{proof}

\begin{prop}
The accumulation points of ${\CCC}_p $ are parameters belonging
 to   $\cup^{p-1}_{k=1}  {\CCC}_k  \cup  \{\infty\}$.
For $p \geq 1$, let ${\lambda}_n\in {\CCC}_p,\, n \in \ZZ,$
be a sequence of virtual centers  with itineraries ${\bf n}_{p}
=(n_1, n_2,\ldots, n_p)$.
\begin{itemize}
\item if the entries $n_i, i=1,2,\ldots,p-1 $ are the same for all
    ${\lambda}_n$ and $n_p=n,$ then the accumulation point of
    ${\lambda}_n$ belongs to ${\CCC}_0=\{\infty\}$ and
\item if the entries $n_i, i=2,\ldots,p  $ are the same for all
    ${\lambda}_n$ and $n_1=n$, then the accumulation point of
    ${\lambda}_n$  is a virtual center ${\lambda} \in {\CCC}_{p-1} $
    with  itinerary ${\bf n}_{p-1}=(n_2,\ldots, n_p)$; that is,
    $f^{p-2}_{\lambda}(\lambda i )=s_{n_2}$.
\end{itemize}

\end{prop}
\begin{proof} Consider the map
$g:\CC  - \cup^{p-1}_{k=1}{\CCC}_k  \to \hat{\CC} $ defined
by $g({\lambda})= f^p_{\lambda}(\lambda i).$ The map   $g$
is not defined at the points $\lambda  \in
\cup^{p-1}_{k=1}{\CCC}_k$ since  these points are essential
 singularities of $g$; $g$ has poles at
$\lambda \in  {\CCC}_{p}$ and is holomorphic otherwise.
 Let $\lambda={\lim}_{n \to \infty}(\lambda_n), \,
 \lambda _n \in {\CCC}_{p} $.
If $\lambda  \notin  \cup^{p-1}_{k=1}{\CCC}_k \cup  \{\infty\}$
 then $g(\lambda)$   is well-defined  and holomorphic  in a
neighborhood of  $\lambda$. On the other hand $\lambda$ is
an accumulation point of  poles and so $g$ has a non-removable
singularity at $\lambda$. Thus we arrive at contradiction.

The proof of the second part of the proposition is the same as
the proof of proposition ~\ref{thm:poleitn}, since  the itineraries of the
virtual centers are
defined by the itineraries  of the corresponding prepoles.
  \end{proof}

\begin{cor}
 If $ \lambda \in \CC -  {\CCC}$
 then for each $p \in \NN $, there is a neighborhood
$U$ of $\lambda$ such that
$$U \cap \left(\cup^{p-1}_{k=1}  {\CCC}_k\right)=  \emptyset.$$
\end{cor}

\section{Deformation of cycles}

We have seen that the periodic points and their multipliers are holomorphic
functions of the parameter $\lambda$, defined by a holomorphic motion, in any
hyperbolic component of the 
parameter space.  We could also obtain these functions by analytic
continuation of local  solutions to the  
 functional equations $\fl{p}{\lambda}(z) = z$ for all positive
integers $p$.  At a boundary point of the hyperbolic component, some
of these local solutions may still have analytic continuations while
others do not.  

Let $\lambda(t), \, t \in [0,1]$ be a path defined in the parameter plane 
$X$ with
$\lambda(0)=\lambda_0$  and $\lambda(1)=\lambda^* \neq \infty$.  For a
given cycle, $\{z_i\}_{i=0}^{p-1}$,  assume that we have a path in $X
\times \C$ defined by analytic continuation of function elements of the
functional equation, 
$z_i(t) =z_i(\lambda(t)) $, $t \in [0,1)$.  Denote the  multiplier function
for the cycle by $m(t)=m(\lambda(t))$. 

 If $\lim_{t \to 1}m(t)=1$,  the functions
$z_i(t)$ cannot be analytically continued to $\lambda^*$ and $\lambda^*$ is
called an {\em algebraic singularity} in analogy with the rational map
case.  Note that this limit is independent of the path.  If, however,   for some path, $\lambda(t) \to \lambda^*$, $\lim_{t \to 1}m(t)= \infty$ or doesn't exist, then $\lambda^*$ is
called a {\em transcendental singularity} for the cycle $z(t)$.

Next we show that a transcendental singularity $\lambda^*$ of a periodic
cycle of period $p$ is a special point  on the boundary of a hyperbolic
component;  it is, in 
fact,  the virtual center of a pair hyperbolic pair, $\Omega_p,\Omega_p'$.  
This, together with theorem~\ref{thm:centers} shows that there is a 
 one to one correspondence between the 
transcendental 
singularities of the periodic points and the virtual centers of the hyperbolic
component pairs.

\begin{lemma}
\label{thm:limitpt} Suppose  $\lambda^*$ is a transcendental
singularity for the cycle $z(t)$.   Then, 
each of the curves $z_i(t)$, $i=0,\ldots,p-1, \, t \in [0,1)$ has a limit
point, perhaps infinite, in $(\lambda^*,\Chat)$, the dynamic plane of
$f_{\lambda^*}$. 
\end{lemma} 

\begin{proof}  Let $A$ be the accumulation set of $\{z_{p-1}(t)\}$ in
$(\lambda^*,\Chat)$.   If $A$ is not a
point, then $A$ is a continuum.  We claim it cannot be a continuum.

Let $w \neq \infty \in A$ and let $t_n \to 1$ be such that
$ z_{p-1}(t_n) \to w$.   Because of the functional relation, 
 $z_{p-1}(t_n)=
\f{p}(z_{p-1}(t_n))$.  As we let $n \to \infty$, either we have
$w=\fl{p}{\lambda^*}(w)$ and $\lambda^*$ is not a transcendental
singularity or $\fl{p}{\lambda^*}(w)$ is not defined and $w$ is a prepole
or $f_{\lambda^*}$ of order at most $p$.  Since there are only countably
many prepoles of order at most $p$, they form a discrete set and $A$ is not
a continuum. The same is true for each path $z_i(t),\, i=1
\ldots, p-1$.  Assume therefore that the endpoint $z_0(1)$ is a prepole so
that  some $z_k(t)$ is
unbounded.  We shall see below that only one curve can be unbounded if the
cycle is not symmetric and that two curves are unbounded if it is
symmetric.  For definiteness, we shall always assume the points are ordered
so that 
$z_{p-1}(t)$ is unbounded.
 \end{proof}

\begin{prop}
\label{thm:singularities}
If the periodic point $z(\lambda)$ has a transcendental singularity
 at $\lambda^*$ then there exists  a
   hyperbolic component $\Omega_p$ (or $\Omega_p'$) with virtual center at
$\lambda^* \in \CCC_{p-1}$  and $\lambda(t)\in \Omega_p$ (resp.
$\Omega_p'$) for $t> t_0$.
\end{prop}

\begin{proof}  Since we will only be interested in the situation where $p$ is
large, we omit the discussion of the special but not too difficult
exceptional cases  $p \leq 2$ and assume $p > 2$. 

By assumption there is a path $\lambda(t)\to \lambda^*, t \in [0,1]$ in
the parameter plane such that the branch $z_{p-1}(t) \to 
\infty$ as $t \to 1$. 
  First we prove that the limit $ z_{p-2}(1)$ is some
pole $s_n $.  We know $$\lim_{t\to 1} z_{p-1}(t)=\lim_{t\to 1} \lambda(t)\tan
z_{p-2}(t)= \infty,$$  $\lambda(t), \, t \in [0,1]$ is bounded and
$\lim z_{p-2}(t)$ is a point in $\CC$.  This implies that $$ \lim_{t\to 1}\tan
z_{p-2}(t)= \infty$$  and that the endpoint of the
curve is a pole $s_n $; that is, $$\lim_{t\to 1} z_{p-2}(t)=\lim_{t\to 1}
f^{-1}_{n,t}( 
z_{p-1}(t))= s_n .$$ 

  Note that 
 the closed curve $z_{p-2}(t) $ is bounded and has ends $z_{p-2}(0)$ and 
$z_{p-2}(1)=s_n$.  It cannot intersect the
curve $\lambda(t) i$ except at the endpoint $\lambda^* i$ because  all
points in $z_{p-2}(t)$ are images of points in $z_{p-3}(t)$ and $\lambda(t)
i$ is omitted for $f_t$.   

We now show that the endpoint of one of the branches is an asymptotic
value.  To do this we need to find an asymptotic path in the dynamic plane
of $f_{\lambda^*}$.  
Draw the projections of the branch curves $z_i(t)$ to this plane and use
the same notation for them.  Now let $\beta(t) =
\fl{p}{\lambda^*}(z_{p-2}(t))$.  Initially, $\beta$ and $z_{p-2}$ may be
far apart, but since $f_{\lambda}(z)$ is holomorphic in both $\lambda$ and
$z$, and since  both paths end at the pole $s_n$, given $\eta$ there is a
$t_0$ such that for $t > t_0$, 
\begin{eqnarray}
 \label{eqn:limalpha} 
|\beta(t) - z_{p-2}(t)|
=|\fl{p}{\lambda^*}(z_{p-2}(t))- \fl{p}{\lambda}(z_{p-2}(t))|
 < \eta.
\end{eqnarray}

Now set $\alpha(t)=f_{\lambda^*}(\beta(t))$.  
Since $\lim_{t\to 1}\beta(t)$ is a pole,   $\alpha(t)$ is unbounded.  The
curve $z_{p-2}(t)$ 
is bounded and we shall assume for the moment for simplicity, that all the
curves $z_i(t)$, $i=0,\ldots,p-2$ are bounded --- we shall eliminate this
assumption at the end of the proof.  Because the preimage of a closed
bounded set under a tangent map is unbounded if and only if the bounded set
contains an asymptotic value, $f_{\lambda^*}(z_{p-1}(t))$
ends at an asymptotic value and thus so does $z_{0}(t)$; furthermore, both
$\alpha(t)$ and $z_{p-1}(t)$ are asymptotic paths; to wit, $\Im z_{p-1}(t)
\to \infty$.  Thus, $\fl{p-2}{\lambda^*}(\lambda^* i)= \pm s_n$ and by
symmetry both $\lambda^*$ and $-\lambda^*$ are virtual centers.

Using the map $f_{\lambda^*}$ and the itinerary of the cycle, the
$i$-th preimage of the curve $\alpha(t)$ has the same limit as $t \to
1$ as the curve $z_{p-1-i}(t)$, $i=1,\ldots, p-1$.

Next we need to prove that $\Re z_{p-1}(t)$ stays bounded.
  The argument is reminiscent of
the proof that there are no Herman rings for these functions.  It isn't
the same, however, because now parameter values vary.  We therefore use
the same trick as we did above;  we use the curve $\alpha(t)$ and work
in the dynamic plane of $f_{\lambda^*}$.  Then, 
if we assume  $\Re z_{p-1}(t)$ is unbounded, $\Re\alpha(t)$ is also
unbounded.  It follows that the range of $\arg(f_{\lambda^*}(\alpha(t)) -
\lambda^* i)$ contains either a semi-infinite interval 
$(\theta,\infty)$ or a union of semi-intervals
$(\theta_{k},\theta_{k+1}),\,k\to \infty$.  Let $U$ be a neighborhood of
$\lambda^* i$ and let $I$ be any straight line through $\lambda^* i$. Then
$f_{\lambda^*}(\alpha(t)) $ meets $I$ infinitely often on both sides of
$\lambda^* i$.

  Now for $r \gg 0$,  $J = f_{n,\lambda^*}^{-1}((-\infty,r)\cup(r,\infty))$ 
 is a curve containing $s_n$ in its interior.  Since
$f_{\lambda^*}^{p-2}(U)$ contains a neighborhood of $s_n$, the curve
$f_{\bf n_{p-2},\lambda^*}^{-(p-2)}(J)$ intersects $f_{\lambda^*}(\alpha(t)$ infinitely often on both
sides of $\lambda^* i$ and hence there is a sequence $t_n$, $n \to \infty$
such that $f_{\lambda^*}^{p-1}(\alpha(t_n)) \cap J \neq \emptyset$ and
$f_{\lambda^*}^{p}(\alpha(t_n)) \in {\RR}$.  This contradicts
equation~(\ref{eqn:limalpha}) and the fact that $\Im z_{p-1}(t)\to \infty $.
We conclude not only that $\Re z_{p-1}(t)$ stays bounded, but that for $t$
close to $1$, it must stay inside a strip of width $\pi$.

To show that the path $\lambda(t)$ enters a hyperbolic component we prove 
 $$
     \lim_{t\to 1} \left|\prod^{p-1}_{i=0}f^{\prime}_t(z_i(t))\right|=0
$$
   where, as usual,  $f_t=f_{\lambda(t)}$.
 We have
\begin{eqnarray}
\label{eqn:mu20}
  |f^{\prime}_t(z_{p-2}(t))|=\left|\lambda(t)\left(
 1+{z^2_{p-1}(t)}/{\lambda^2(t)}\right)\right|
\approx {\rm O}( |\lambda(t)|^{-1}|z^2_{p-1}(t)|)
\end{eqnarray}
and
\begin{eqnarray}
\label{eqn:mu21}
  \left| f^{\prime}_t(z_{p-1}(t))\right|=|{\lambda}(t)|
   |\sec^{2}( z_{p-1}(t))|= 4 |{\lambda}(t)|\left|
   (e^{i z_{p-1}(t)}+ e^{-i z_{p-1}(t)})^{-2}\right|.
\end{eqnarray}
    Since $|\lambda(t)|<K$ for some $K>0$ and  $\Im z_{p-1}(t)\to\infty $
 as $t \rightarrow 1$, 
\begin{eqnarray}
\label{eqn:mul11}
   |\sec^{2}( z_{p-1}(t))|={\rm O}\left( |{\lambda}(t)|
        e^{-2\Im z_{p-1}(t)}\right) \rightarrow 0
\end{eqnarray}

 By equations~(\ref{eqn:mu20})---~(\ref{eqn:mul11}), as $t \to 1$,
$$|\sec^2 z_{p-2}(t)\sec^2 z_{p-1}(t))| =
 {\rm O}(|z_{p-1}(t)|^{2} e^{-2\Im z_{p-1}(t)}).$$
Because  $\Re z_{p-1}(t)$ is universally
 bounded for $t\to 1$ , while $\Im z_{p-1}(t)\to \infty $,

$$ |z_{p-1}^2(t)|e^{-2\Im z_{p-1}(t)}\to 0.$$
By assumption,  the curves $z_i(t), \, i=0, \ldots, p-2$
  stay  a
positive distance from all poles (at least for $t$ close to 1), so
 there exists a bound $C>0 $ such that
$$
      \left|\prod^{p-3}_{i=0}f^{\prime}_t(z_i(t))\right|\leq C.
$$
Then
$$ |m(t)|={\rm O}( \Im z^2_{p-1}(t)e^{-2\Im z_{0}(t)})\to 0$$
so there is some $t' \in [0,1)$ such that the  multiplier
  $$ |m(t)|=\left|\prod^{p-1}_{i=0}f^{\prime}_t(z_i(t))\right|<1 $$
 for $t>t'$; $ z_i(t'), i=0,\ldots,p-1$ is a neutral
 periodic cycle. Thus for $ t>t',  \, \lambda(t)  $   belongs to
 some hyperbolic component ${\Omega}_p$ of $\HHH^0$ and $\lambda^* =
 \lim_{t\to 1}\lambda(t)$ is a virtual center of $ {\Omega}_p$.

Now suppose there is a branch $z_i(t)\neq z_{p-1}(t)$ such that $\lim_{t
\to 1}z_i(t)=\infty$.  Arguing as above, $z_{i-1}$ tends to a pole and
$z_{i+1}$ tends to an asymptotic value.  Since there are only two
asymptotic values, either the cycle is symmetric and contains both 
asymptotic values or it has a lower period $q|p$.  Again since $\Re z_i(t)$
is bounded, $\sec^2 z_{i-1}(t)\sec^2 z_i(t)\to 0$ and for all branches
$z_j(t)$ such that $z_j(1)$ is not a pole or $ z_j(1)\neq\infty$, there is
a constant independent of $t$ and $j$ such that $|\sec^2z_j(t)|<C$; thus, as
$t\to 1$ the periodic cycle becomes an attractive cycle of the same period.

Finally, from the proof of theorem~\ref{thm:centers}, we conclude either
$f_t$ has a distinct symmetric cycle $-z_i(t), i=0, \ldots, p-1$
and that $q=p$, or the cycle $z_i(t), i=0,\dots, p-1$ is symmetric with
respect to the origin,   $q=p/2$ and $\lambda\in \Omega_{p'}$ with
virtual center at $\lambda^*$.   
\end{proof}

We close this discussion by considering the  behavior of a periodic
cycle  and its multiplier,  in and on the boundary of the hyperbolic
component where it is the attracting cycle.  

\begin{prop}
\label{thm:newprop} Let $\lambda^*$ be the virtual center of the
component pair $(\Omega_{p},\Omega_{p}').$  If  $z_{0}(\lambda),
z_1(\lambda), \ldots, z_{p-1}(\lambda)$ is
an attracting periodic cycle of period $p$
 for $\lambda \in \Omega_p$, and if  $m_{\lambda}$ is its
multiplier,
then
\begin{enumerate}
\item $z_i(\lambda)$ and $m_{\lambda}$ are  transcendental
 meromorphic functions of $\lambda$;\\
 $m_{\lambda}$ has  an essential
 singularity at $\lambda^*$ as does  $z_{\lambda} = z_k(\lambda)$ for some
$k \in \{ 0,\ldots,p-1\}$.

\item $\lambda^* i$ is an asymptotic value of
$f_{\lambda^*}=z_{k+1}(\lambda^*)$,
$ 0$ is an asymptotic value of $m_{\lambda}$  and $\Omega_p$
is an asymptotic tract for both $m_{\lambda}$ and $z_{k+1}(\lambda)$.

\item $z_{\lambda}$ has algebraic singularities at a sequence of
 points $\lambda_j\in \partial{\Omega}_p $ where
 $m(z_{\lambda_j})=1$ and $\lambda_j \to \lambda^*$ as
 $j \to \infty$.
\end{enumerate}
\end{prop}

\begin{proof} By proposition~\ref{thm:singularities}, if $\lambda \to \lambda^*$
inside $\Omega_p$, then for some $k$, $\Im
z_k(\lambda) \rightarrow \infty$ and $f_{\lambda}(z_k(\lambda)) \rightarrow
\lambda^* i$. Relabelling the cycle if necessary we may assume $k=p-1$ so
that $z_{0}(\lambda) \rightarrow \lambda^* i$.
Now set
 $z_{\lambda} =z_{0}(\lambda)$.
  Because $\lambda^* i$ is
a prepole of $f_{{\lambda}^*}$ of order $p-1$,
$\fl{p}{\lambda^*}(\lambda^* i)$ is not defined and so neither is
$m_{\lambda}$, proving 1.

To prove 2. note that
the equation $\fl{p}{\lambda}(z_{\lambda})=z_{\lambda}$ has no
solution for $\lambda=\lambda^*$ and that ${\lambda}^* i$ is
an asymptotic value of $z_{\lambda}$.
Because
$\lambda^*$ is a
virtual center of ${\Omega}_p,$ we see from proposition ~\ref{thm:prop2.14}
 that if
$\lambda  \to {\lambda^*}$ inside ${\Omega}_p$ then $m_{\lambda}
\rightarrow 0$.

The third point
follows from propositions ~\ref{thm:prop2.14} and ~\ref{thm:singularities}.
\end{proof}

\bibliography{ref}

\end{document}